\documentclass[12pt,twoside]{article}
\usepackage[latin1]{inputenc}
\usepackage{times}
\usepackage{epsfig}
\usepackage{dsfont}
\usepackage{amsfonts}
\usepackage{float}
\usepackage{amsmath}
\usepackage{amsthm}
\usepackage{latexsym}
\usepackage{xr} 
\usepackage{natbib}
\usepackage{xcolor}
\usepackage{multirow}
\usepackage{rotating}
\usepackage{enumitem}

\relax 
\newcommand{\pkg}[1]{{\fontseries{b}\selectfont #1}}

\allowdisplaybreaks[3]

\newcommand{\1}{{\mathds{1}}}

\renewcommand{\qed}{\mbox{ } \hfill $\Box$\\ }
\newcommand{\mc}{\multicolumn}

\newcommand{\bay}{\begin{array}}
\newcommand{\eay}{\end{array}}

\newcommand{\bqa}{\begin{eqnarray*}}
\newcommand{\eqa}{\end{eqnarray*}}

\newcommand{\bee}{\begin{eqnarray*}}
\newcommand{\eee}{\end{eqnarray*}}

\newcommand{\bea}{\begin{eqnarray*}}
\newcommand{\eea}{\end{eqnarray*}}

\newcommand{\bqan}{\begin{eqnarray}}
\newcommand{\eqan}{\end{eqnarray}}

\newcommand{\be}{\begin{eqnarray}}
\newcommand{\ee}{\end{eqnarray}}

\newcommand{\bit}{\begin{itemize}}
\newcommand{\eit}{\end{itemize}}

\newcommand{\ben}{\begin{enumerate}}
\newcommand{\een}{\end{enumerate}}

\newcommand{\beq}{\begin{equation}}
\newcommand{\eeq}{\end{equation}}

\newcommand{\bdes}{\begin{description}}
\newcommand{\edes}{\end{description}}

\newcommand{\btb}{\begin{tabular}}
\newcommand{\etb}{\end{tabular}}

\newcommand{\bcen}{\begin{center}}
\newcommand{\ecen}{\end{center}}

\newcommand{\bmp}{\begin{minipage}}
\newcommand{\emp}{\end{minipage}}

\newtheorem{thm}{Theorem}

\newtheorem{cor}{Corollary}

\newtheorem{cond}{Condition}

\newcommand{\rnc}{\renewcommand}
\newcommand{\nc}{\newcommand}
\newcommand{\mrm}{\mathrm}
\renewcommand{\b}{\textbf}

\newcommand{\bs}{\boldsymbol}
\nc{\mb}{\mathbb}
\nc{\mac}{\mathcal}
\nc{\E}{\mb{E}}
\nc{\N}{\mb{N}}
\nc{\R}{\mb{R}}
\nc{\Q}{\mb{Q}}
\rnc{\P}{\mrm P}
\nc{\bP}{\b P}
\nc{\bA}{\b A}
\rnc{\d}{\mrm d}
\nc{\C}{\mc{C}}
\nc{\D}{\mc{D}}
\nc{\B}{\mc{B}}
\nc{\I}{\mc I}
\nc{\J}{\mc J}
\nc{\JI}{{\J\I}}
\rnc{\IJ}{{\I\J}}
\nc{\ji}{{\J|\I}}
\nc{\gDg}{\stackrel{d}{=}}
\nc{\oPo}{\stackrel{\mrm p}{\longrightarrow}}
\nc{\oWo}{\stackrel{w}{\longrightarrow}}
\nc{\oDo}{\stackrel{d}{\longrightarrow}}
\nc{\nae}{Nelson-Aalen estimator}
\nc{\aje}{Aalen-Johansen estimator}
\nc{\naeL}{Nelson-Aalen estimator\ }
\nc{\ajeL}{Aalen-Johansen estimator\ }
\nc{\CIF}{cumulative incidence function}
\nc{\CIFL}{cumulative incidence function\ }
\nc{\wh}{\widehat}

\newcommand{\kp}{\otimes}

\newcommand{\Cov}{\operatorname{{\it Cov}}}
\newcommand{\Var}{\operatorname{{\it Var}}}

\newcommand{\vnull}{{\bf 0}}

\newcommand{\pu}{\color{black}}
\newcommand{\bl}{\color{black}}
\definecolor{darkolivegreen}{rgb}{0.33, 0.42, 0.18}
\definecolor{darkgoldenrod}{rgb}{0.72, 0.53, 0.04}

\usepackage{lipsum}

\newcommand{\g }{\color{black}}

\begin{document}


\title{\Large \bf Nonparametric MANOVA in Mann-Whitney effects$^{1}$
}
\author{Dennis Dobler$^{*}$, Sarah Friedrich$^\dagger$ and  Markus Pauly$^\dagger$ \\[1ex] 
}
\maketitle

\begin{abstract}
{\color{black}
\noindent Multivariate analysis of variance (MANOVA) is a powerful and versatile method to infer and quantify main and interaction effects in metric multivariate multi-factor data. It is, however, neither robust against change in units nor a meaningful tool for ordinal data. Thus, we propose a novel nonparametric MANOVA. Contrary to existing rank-based procedures we infer hypotheses formulated in terms of meaningful Mann-Whitney-type effects in lieu of distribution functions. The tests are based on a quadratic form in multivariate rank effect estimators and critical values are obtained by the bootstrap. This newly developed procedure provides asymptotically exact and consistent inference for general models such as the nonparametric Behrens-Fisher problem as well as multivariate one-, two-, and higher-way crossed layouts.}
{\pu Computer simulations in small samples confirm the reliability of the developed method {\color{black} for ordinal as well as metric data with covariance heterogeneity}. Finally, an analysis of a real data example illustrates the applicability and correct interpretation of the results.
}
\end{abstract}

{\pu Multiple Samples;}
Rank-based Methods;
{\g Wild} Bootstrap.

\vfill
\vfill

\noindent${}^*$ Vrije Universiteit Amsterdam, Department of Mathematics, Netherlands\\
 {\pu \mbox{ }\hspace{1 ex}email: d.dobler@vu.nl} 
\medskip \\
\noindent${}^{\dagger}$ {Ulm University, Institute of Statistics, Germany\\
 {\pu \mbox{ }\hspace{1 ex}email: sarah.friedrich@uni-ulm.de}\\
 \mbox{ }\hspace{1 ex}email: markus.pauly@uni-ulm.de}

\noindent${}^{1}$ Authors are in alphabetical order.

\newpage


\section{Motivation and Introduction}\label{int}

In many experiments, e.g., in the life sciences or in econometrics, observations are obtained in elaborate factorial designs with multiple endpoints. Such data are usually analyzed using MANOVA methods such as Wilk's $\Lambda$. These procedures, however, rely on {\pu the assumptions of} multivariate normality  and covariance homogeneity and usually break down if these {\pu prerequisites} are not fulfilled. 
In particular, if the observations are not even metric,  such {\bl applications} are no longer possible
since means no longer provide adequate effect measures. 
{\g To this end, several rank-based methods have been proposed for nonparametric MANOVA and {\pu repeated measures designs} which are usually based on Mann-Whitney-type effects:}
In the context of a nonparametric {\it univariate} two-sample problem with independent {\color{black}and continuous} observations $Y_{ik} \sim F_i, \ i=1, 2, \ k=1, \dots, n_i$, \citet{mann1947test} introduced the effect $w = P(Y_{11} \leq Y_{21})=\int F_1 dF_2$ {\bl also known as ordinal effect size measure \citep{Acion2006}.} 
An estimator of $w$ is easily obtained by replacing the distribution functions with their empirical counterparts. 
While this effect has several desirable properties and is widely accepted in practice \citep{brumback2006using, Kieser2013}, 
generalizations to more than one dimension or higher-way factorial designs are not straightforward.

{\color{black} Concerning the latter, }{\pu there {\g basically} exist two possibilities in the literature to cope with $ a \geq 3 $ sample groups with independent {\color{black} univariate} observations $Y_{ik} \sim F_i, \ i=1, \dots, a, \ k=1, \dots, n_i$:
	{\it First}, considering only the pairwise effects $w_{i\ell} =P(Y_{i 1} \leq Y_{\ell 1}), \ {\pu  1 \leq i \neq \ell \leq  a }$ 
	(as proposed by \citealp{rust1984modification}) can lead to paradox results in the sense of Efron's Dice;
	see also \cite{thas12} and the contributed discussions by M. P. Fay and W. Bergsma and colleagues 
	for pros and cons of the possibly induced intransitivity  by {\color{black} certain} probabilistic index models. 
	We refer to \citet{brown2002,THANGAVELU2007} or \citet{brunner2016JRSSB} and the references cited therein for {\color{black} further} considerations 
	{\bl on this issue.}
{\g \it Second,} in order to circumvent the problem of {\pu intransitive} effects, the {\color{black} group-wise} distribution functions {\color{black} $F_i$ may} be compared to the same reference distribution.
{\pu Usually, this is} the pooled distribution function $H =\frac{1}{N}\sum_{i=1}^a n_i F_i$ \citep{kruskal1952nonparametric,kruskal1952use}, resulting in 
{\bl so-called \citep[e.g.,][]{brunner2016JRSSB}} relative effects $r_i = \int H dF_i$. 
{\color{black} Multivariate generalizations of this approach can be found in \citet{puri1971nonparametric}, \citet{munzel2000nonparametric} or \citet{BrunnerMunzelPuri2002};}
{\color{black}see also \cite{deneve15} for a related approach.}
Since these quantities depend on the sample sizes $n_i$, however, they are no fixed model constants and changing the sample sizes {\pu might} dramatically {\pu alter} the results; see {\pu again} \citet{brunner2016JRSSB} for an example {\bl in the univariate case}. 
{\pu For this reason}, \citet{brunnerPuri} proposed a different {\bl nonparametric }effect $p_i =\int G dF_i$ {\color{black} for univariate factorial designs,}
where $G=\frac{1}{a}\sum_{i=1}^a F_i$ denotes the unweighted mean of {\pu all} distribution functions. 
{\bl The same approach has also been extended to other settings by \cite{gao05}, {\color{black}\citet{gao08-2}}, \cite{gao08} 
and  \cite{UKP2016}.} Nevertheless, none of them considered null hypotheses formulated in terms of fixed and meaningful model parameters.
For a more intuitive interpretation of the results, however, it is sensible to formulate and test hypotheses in more vivid {\bl effect sizes.} 
In particular, it is widely accepted in quantitative research that ``\emph{effect sizes are the most important outcome of empirical studies}'' \citep{lakens13}.
{\color{black}\citet{brunner2016JRSSB}} therefore infer null hypotheses {\pu stated} in terms of the unweighted {\bl nonparametric} effects {\pu via} $H_0^p: \b H \b p =\vnull$ for a suitable hypothesis matrix $\b H$ and the pooled vector $\b p$ of {\pu the effects $p_i$}{\color{black}, see also \citet{konietschke2012rank} for the 
special case of one group repeated measures.}
}

{\pu 

	In the present paper, we strive to generalize their models and methods in several directions:
	\begin{enumerate}
		\item We examine generalizations to the more involved context of multivariate data where dependencies between observations from the same unit {\pu need} to be taken into account.
		This multivariate case allows for testing hypotheses on the influence of several factors on single or several outcome measurements.
		\item {\bl More general as in Repeated Measures designs the} outcomes in different components may be measured on different units (such as grams and meters){\bl. In particular, they actually} need not even be elements of metric spaces; 	totally ordered sets serve equally well as spaces of outcomes because we develop rank-based methods for our analyses.
		Query scores are an example of such ordered data without having a unit in general.
		Differences will be tested with the help of {\color{black} a quadratic form in the rank-based effect estimates.}
		\item 
		This test statistic is analyzed by means of modern empirical process theory (instead of the more classical and sometimes cumbersome projection-based approaches for rank statistics).
		Since it is asymptotically non-pivotal, appropriate bootstrap methods for asymptotically reproducing {\color{black} its} correct limit null distribution {\color{black} are proposed}. As bootstrapping entails several good properties when applied to empirical distribution functions
		and our rank-based estimates offer a representation as a functional of multiple empirical distribution functions,
		we {\g expect} to obtain reliable inference methods using bootstrap {\color{black} techniques}.
		This conjecture will be supported by simulation results which indicate a good control of the type-I error rate even for small sample set-ups {\color{black} with ordinal or heteroscedastic metric data}.
	\end{enumerate}
}
{\g Our model formulation thereby comprises novel procedures for general multivariate factorial designs with crossed or nested factors and even contains the so-called nonparametric multivariate Behrens-Fisher problem as a special case.} {\bl Moreover, the methodology also allows for subsequent post-hoc tests}.

The paper is organized as follows: In Section \ref{mod} we describe the statistical model and the null hypotheses of interest. Section \ref{sec:asy} presents the asymptotic properties of our estimator
{\pu and, subsequently, states the asymptotic validity of} {\bl its bootstrap versions. Deduced statistical inference procedures are discussed in Section~\ref{inf} and their small sample behavior is} analyzed in extensive simulation studies in Section \ref{sim}. Section \ref{app} contains the {\pu real data} analysis of {\pu the gender influence on education and annual household income of shopping mall customers in the San Francisco Bay Area}. 
We conclude with some final remarks in Section \ref{dis}.
{\pu The proofs of all theoretic results and the derivation of the asymptotic covariance matrices are given in the Appendices~\ref{app:proofs} and~\ref{app:cov}, respectively.}
Proof of all large sample properties of the classical bootstrap applied in the present framework, further simulation results regarding the power of all proposed methods, and additional analyses of the shopping mall customers data example are provided in Appendices~\ref{sec:cbs}, \ref{sec:power}, and \ref{sec:add_ana}, respectively.

\newpage

\section{Statistical Model} \label{mod}

{\pu Throughout, let $(\Omega, \mac A, P)$ be a probability space on which all random variables will be defined.}
We assume a general factorial design with multivariate data, that is, we consider independent random vectors
\bqan \label{eq:model}
\b X_{ik} =(X_{ijk})_{j=1}^d {\pu : ~ \Omega \longrightarrow \R^d}, ~ i=1, \dots, a; ~ k = 1, \dots, n_i
\eqan
{\color{black}of dimension $d\in\N$,} where $X_{ijk}$ denotes the $j$-th measurement of individual $k$ {\pu in} group $i$. 
{\pu Thus, the total sample size is $N = \sum_{i=1}^a n_i$. 
The distribution of $\b X_{ik}$ is assumed to be the same within each group with marginals denoted by}
\[
X_{ijk} \sim F_{ij}, \hspace{0.5cm} i=1, \dots, a, ~k=1, \dots, n_i, ~j=1, \dots, d.
\]


{\pu Throughout, we understand all $F_{ij}$ as the} so-called normalized distribution functions, {\pu i.e. the means of their left- and right-continuous versions} \citep{ruymgaart1980, akritas1997nonparametric, munzel1999}. 
{\bl This allows for a unified treatment of metric and ordinal data and will later on lead to statistics formulated in terms of mid-ranks}. 
For convenience, we {\pu combine} the observations $\b X_{ik}$ in larger vectors
\bqan
\b X_i \  = \ (\b X_{i1}', \dots, \b X_{in_i}')' {\pu : ~ \Omega \longrightarrow \R^{d n_i} }, \quad \text{and}\quad \b X &= & (\b X_1', \dots, \b X_a')' {\pu : ~ \Omega \longrightarrow \R^{d N} }
\eqan
containing all the information of group {\pu $i=1,\dots,a$} and the pooled sample, respectively. 
{\bl Different to the special case of repeated measurements \citep{konietschke2012rank, brunner2016JRSSB}} the components are in general not commensurate. Therefore, 
comparisons between the different groups {\pu are} performed component-wise. 
{\pu To this end,} let 
$G_j = \frac{1}{a} \sum_{i=1}^a F_{ij}, \, j=1, \dots, d$ denote the unweighted mean distribution function for the $j$-th component. 
We consider $G_j$ as a benchmark distribution for {\pu comparisons in} the $j$-th component. 
{\pu In particular,} denote by $Y_j\sim G_j$ a random variable 
that is independent of $\b X$ and define unweighted {\bl nonparametric} effects for group $i$ and component $j$ by 
\bqan \label{trteffects}
p_{ij} = P(Y_j<X_{ij1}) + \tfrac{1}{2}P(Y_j=X_{ij1}) = \int G_j dF_{ij} = {\pu \frac1a \sum_{\ell = 1}^a w_{\ell ij} = ~ } \overline{w}_{\cdot ij},
\eqan
where $w_{\ell ij} = \int F_{\ell j} dF_{ij} {\color{black} =P(X_{\ell j 1} < X_{i j1}) + \frac{1}{2} P(X_{\ell j1} = X_{ij1})}$ {\pu quantifies the {\color{black} Mann-Whitney} effect for groups $\ell$ and $ i$ in component $j$.
Note that $w_{\ell ij} = 1/2$ in case of $\ell = i$.}
This definition naturally extends the univariate effect measure given in \cite{brunner2016JRSSB} to our {\color{black} general} multivariate {\pu set-up}. 
Note that, {\pu in contrast to their suggestion for an extension to repeated measures designs,} 
comparisons with respect to the overall mean distribution $G = \frac{1}{ad} \sum_{i=1}^a\sum_{j=1}^d F_{ij}$ {\pu are not appropriate here 
since  we {\bl study a more general model that allows} for components measured on different units. 
{\pu However,} the advantages of an unweighted effect measure as discussed in \citet{brunner2016JRSSB} {\pu still apply: 
The {\color{black} $p_{ij}$'s in \eqref{trteffects}} are fixed model quantities that do} not depend on the sample sizes $n_1, \dots, n_a$,
{\pu thus allowing} for {\color{black}a }transitive ordering. 
Moreover, interpretation {\g of these effects} is rather simple: An effect $p_{ij}$ smaller than 1/2 means that observations from the distribution $F_{ij}$ (i.e. 
from component $j$ in group $i$) tend to smaller values than those from the corresponding {\it benchmark distribution} $G_j$. \\
In this set-up, we formulate null hypotheses as $H_0^{p} : \b H \b p = \bs {0}$
where $\b p = (p_{11}, {\pu p_{12},} \dots, p_{ad})'$ denotes the vector of the relative effects $p_{ij}$, $i = 1, \dots, a, ~ j=1, \dots, d$ and $\b H$ is a suitable hypothesis matrix {\pu with $ad$ columns}. Instead of $\b H$ we may {\pu equivalently} use the unique projection matrix $\b T = \b H'(\b H \b H')^+ \b H$
which is idempotent and symmetric and fulfills $\b H \b p = \bs {0} \Leftrightarrow \b T \b p = \bs {0}$; see e.g., \citet{BrDeMu:1997,brunnerPuri} and \citet{brunner2016JRSSB}.
Henceforth, let $\b I_d$ and $\b J_d$ denote the $d$-dimensional unit matrix and the $d\times d$ matrix of 1's, respectively, and define by $\b P_d = \b I_d - \frac{1}{d} \b J_d$ the so-called $d$-dimensional centering matrix.

{\bl In particlar, in case of $a=2$} our approach includes the nonparametric multivariate Behrens-Fisher problem {\bl
$$
H_0^p(\b T) : \{\b T \b p = \b 0 \} = \{ \b p_1 =\b p_2 = \b 1_d / 2 \} 
$$
with $\b T = \b P_2{\bl \otimes \b I_d} =  \frac12 (\begin{smallmatrix} 1 & -1 \\ -1 & 1 \end{smallmatrix}){\bl \otimes \b I_d}$ 
and $\b p_i =(p_{i1},\dots,p_{id})'$, $i=1,2$. Similarly, one-way layouts are covered by choosing $\b T = {\bl \b P_a \otimes \b I_d}$, leading to the null hypothesis 
$H_0^p(\b T): \{\b p_1 = \dots = \b p_a\}$.  Moreover,} more complex factorial designs {\bl can be treated as well } by splitting up the {\pu group} index $i$ into sub-indices $i_1, i_2, \dots$ according to the number of factors considered.
{\g For example, consider a two-way layout with crossed factors $A$ and $B$ with levels $i_1=1, \dots, a$ and $i_2=1, \dots, b$, respectively. 
{\color{black}	In this case, the random vectors in \eqref{eq:model} become $\b X_{i_1 i_2 k}, i_1 = 1, \dots a,~ i_2=1, \dots, b, ~k=1, \dots, n_{i_1i_2}$.
}
	We thus obtain the effect vector $\b p =(\b p_{11}', \dots, \b p_{ab}')'$, where all vectors $\b p_{i_1i_2}{\color{black}=(p_{i_1i_21}, \dots, p_{i_1i_2d})'}, i_1=1, \dots, a, i_2=1, \dots, b$ are $d$-variate and there are $n_{i_1 i_2}>0$ subjects observed at each factor level combination. 
	Hypotheses of interest in this context are the hypotheses of no main effects as well as the hypothesis of no interaction effect between the factors. The hypothesis of no main effect of factor $A$ can be written as $H_0(A): \{(\b P_a \kp \frac{1}{b} \b J_b \kp \b I_d)\b p = \vnull \}$,
	where $\kp$ denotes the Kronecker product.
	Similarly, the hypothesis of no effect of factor $B$ is formulated as  $H_0(B): \{(\frac{1}{a} \b J_a \kp \b P_b \kp \b I_d)\b p = \vnull \}$ and the hypothesis of no interaction effect as  $H_0(AB): \{(\b P_a \kp  \b P_b \kp \b I_d)\b p = \vnull \}$.}
{\color{black} For other covered factorial designs and corresponding contrast matrices we refer to Section 4 in \citet{Kon:2015}.}
Equivalent formulations of the above null hypotheses in terms of the illustrative but notationally more elaborate decomposition into all factor influences are given in the Supplementary Material of \cite{brunner2016JRSSB} {\color{black} for the univariate case, but directly carry over to the present context}.
We note that in the general multivariate case, null hypotheses like $H_0^p$ have only been considered 
in the special case of the nonparametric Behrens-Fisher problem \citep{BrunnerMunzelPuri2002}.
Up to now, {\bl multivariate testing procedures} for one-, two-, or even higher-way layouts focus on null hypotheses formulated in 
terms of distribution functions; see, e.g., \cite{bathke08}, \cite{harrar08}, \cite{harrar12} and the references given in Section~\ref{int}.

To estimate the vector of effects, 
we consider the empirical {\pu (normalized)} distribution functions $\wh F_{ij}(x) = \frac{1}{n_i}\sum_{k=1}^{n_i} c(x - X_{ijk})$ 
where $c(u) = \1\{ u > 0 \} + \frac12 \1\{ u = 0 \}$.
Thus, we obtain estimators for the nonparametric effects $p_{ij}$ by replacing the distribution functions with their empirical counterparts
\bqa
\wh p_{ij} = \int \wh G_j d\wh F_{ij} = \frac{1}{a} \sum_{\ell=1}^a \wh w_{\ell ij},
\eqa
where {\pu $\wh G_j = \frac1a \sum_{\ell=1}^a \wh F_{\ell j}$ and}
$$\wh w_{\ell ij} = \int \wh F_{\ell j} \d \wh F_{ij} = \frac1{n_\ell} \frac1{n_i} \sum_{k=1}^{n_i} \sum_{r=1}^{n_\ell} c(X_{ijk} - X_{\ell j r}) 
= \frac{1}{n_\ell}\left(\overline{R}_{ij\cdot}^{(\ell i)} - \frac{n_i +1}{2}\right). $$

Here, $R_{ijk}^{(\ell i)}$ denotes the (mid-)rank of observation $X_{ijk}$ in dimension $j$ among the $(n_i + n_\ell)$ observations in the pooled sample $X_{\ell j1}, \dots, X_{\ell j n_\ell}, X_{ij1}, \dots, X_{ij n_i}$
and $\overline{R}_{ij\cdot}^{(\ell i)}= \frac{1}{n_i} \sum_{k=1}^{n_i} R_{ijk}^{(\ell i)}$ {\pu are} the corresponding rank means.
{\pu We combine all estimated effect sizes into the $ad$-dimensional vector
$\wh {\b p} = ( \wh p_{11}, \wh  p_{12}, \dots, \wh p_{ad})' $.}
{\bl To detect deviations from n}ull hypotheses {\pu of the form $H_0^p(\b T): \{\b T\b p= \b 0\}$} {\bl we propose the application of the following } ANOVA-type test statistic (ATS) 
\bqan\label{eq:ATS}
T_N = N \widehat{\b p}' \b T \widehat{\b p},
\eqan
where {\pu again} $N = \sum_{i=1}^a n_i$ denotes the total sample size in the experiment.

\section{Asymptotic Properties and Resampling Methods}\label{sec:asy}

 In this section, we discuss asymptotic properties of the {\color{black} vector of {\bl estimated effect sizes} $\widehat{\b p}$ {\pu and {\bl propose bootstrap methods to approximate its unknown limit distribution}. For a lucid presentation of the results we thereby assume the following sample size condition:}
 \begin{cond}
  \label{cond:main}
  $\frac{n_i}{N} \rightarrow \lambda_i \in (0,1)$ for all {\pu groups} $i=1, \dots, a$ as $N \rightarrow \infty$.
 \end{cond}
 In other words, {\bl no group shall constitute a vanishing fraction of the combined sample.
 Due to the Glivenko-Cantelli theorem in combination with the continuous mapping theorem,
 the consistency of $\widehat{\b p}$ for $\b p$ follows already under the  weaker assumption $\min_{1\leq i \leq a}(n_i) \rightarrow \infty$.}
 Asymptotic normality {\bl is established in our main theorem below}:
 \begin{thm}
  \label{thm:clt}
  Suppose Condition~\ref{cond:main} holds.
  As $N \rightarrow \infty$, we have
  \begin{align}
  \label{eq:conv}
   \sqrt{N}(\widehat{\b p} - \b p) \oDo \b Z \sim N_{ad}( \bs 0_{ad}, \bs \Sigma),
  \end{align}
  {\pu where the 
  	{\color{black} rather cumbersome form of the}
  asymptotic covariance matrix $\bs \Sigma \in \R^{ad \times ad}$ is {\color{black} stated} in Appendix~\ref{app:cov}
  and $\bs 0_{ad}  \in \R^{ad}$ denotes the zero vector.}
 \end{thm}
 Note that this theorem immediately implies the asymptotic normality of $\sqrt{N} \b T \widehat{\b p}$ under $H_0^p{\bl(\b T): \{\b T \b p = \b 0\}}$.
 Thus, the continuous mapping theorem {\pu yields} the corresponding convergence in distribution for {\color{black} the quadratic form} $T_N$ {\bl defined in \eqref{eq:ATS}}:
 {\pu 
 \begin{cor}
 \label{cor:ATS}
  Suppose Condition~\ref{cond:main} holds.
  As $N \rightarrow \infty$, we have under $H_0^p{\bl(\b T): \{\b T \b p = \b 0\}}$
  \begin{align}
  \label{eq:conv_TN}
   T_N = N \widehat{\b p}' \b T \widehat{\b p} \oDo  \b Z' \b T \b Z \gDg \sum_{{\g h}=1}^{ad} \nu_h Y_h^2,
  \end{align}
  where $Y_1, \dots, Y_{ad}$ are independent and standard normally distributed 
  and $\nu_1, \dots, \nu_{ad} \geq 0$ are the eigenvalues of $\bs \Sigma^{1/2} \b T \bs \Sigma^{1/2}$.
 \end{cor}
 }
 %
 %
 %
 %
 %
 %

 
%
 
 The limit theorem \eqref{eq:conv} raises the question how to calculate adequate critical values for tests in $T_N$. A first naive idea might be to approximate the right hand side of \eqref{eq:conv} by using its representation as a weighted sum of independent $\chi^2$-variables together with (consistent) estimators for the involved eigenvalues $\nu_h$ or the covariance matrix $\bs \Sigma$. However, these choices usually result in too liberal inference methods as already observed by \citet{brunner2016JRSSB} for the univariate case. {\bl Another idea would be to generalize the $F$-approximation proposed in \citet{brunner2016JRSSB} to the present situation. But since this will in general not lead to asymptotic correct level $\alpha$ tests (even in the most simple univariate two sample setting with $a=2$ and $d=1$, see Brunner et al., 2017), we instead focus on resampling the test statistic $T_N$. In particular, we study} two bootstrap approaches for recovering the unknown limit distribution of $T_N$ under $H_0^p$: A sample-specific as well as a wild bootstrap. 
 	
 Here, {\bl a} wild bootstrap {\bl approach} is implemented in the following fashion: 
 {\bl
 first, we notice that $\sqrt{N} (\widehat{\b p} - \b p)$ has an asymptotically linear representation in $\sqrt{N} ((\wh F_{11}, \dots, \wh F_{ad})' - (F_{11}, \dots, F_{ad})')$.
 Indeed, if we denote $\phi_i(f_1, \dots, f_a) = \int (\frac1a \sum_{\ell=1}^a f_\ell) \d f_i$ for functions $f_1,\dots,f_a$ such that the integral is well-defined,
 then
 \begin{align}
 \begin{split}
 \label{eq:func_delta}
  \sqrt{N} (\widehat{p}_{ij} - p_{ij}) & = \sqrt{N} ( \phi_{i}(\wh F_{1j}, \dots, \wh F_{aj}) - \phi_{i}(F_{1j}, \dots, F_{aj}) )  \\
 & = \sqrt N \int ( \wh G_j - G_j) \d F_{ij} - \sqrt N \int (\wh F_{ij} - F_{ij}) \d G_j + o_p(1).
 \end{split}
 \end{align}
 The second equality follows from the functional delta-method applied to the integral functionals $\phi_i$;
 cf. \cite{dopa17} for the two-sample case.
 In this asymptotic expansion, the now proposed wild bootstrap tries to estimate each involved residual $\varepsilon_{\ell jk}(x) = c(x - X_{\ell jk}) - F_{\ell j}(x), k=1, \dots, n_\ell$
 by another centered quantity with approximately the same conditional variance given $\b X$:
 $$ \widehat \varepsilon_{\ell jk}(x) = D_{\ell k} \cdot [ c(x - X_{\ell jk}) - \widehat{F}_{\ell j}(x) ], \quad k=1, \dots, n_\ell, $$ 
 where $D_{\ell k}, i=1,\dots, a, \ k = 1, \dots, n_\ell,$ are i.i.d.\ zero-mean, unit variance random variables with 
 $\int_0^\infty \sqrt{P(|D_{11}| > x)} \d x < \infty$.
 This condition is implied by $E | D_{11} |^{2 + \eta} < \infty$ for any $\eta >0$, and it is thus a weak assumption on the heaviness of tails; cf.\ p.\ 177 in \cite{vaart96}.
 Note that, for each $\ell,k$, our wild bootstrap implementation uses the same multiplier $D_{\ell k}$ for every component $j$ in order to ensure an appropriate dependence structure.
 
 Additionally, apart from estimating the residuals $ \varepsilon_{\ell jk}$ by $ \widehat \varepsilon_{\ell jk}$,
 the unknown distribution functions $F_{\ell j}$ in the integrators in~\eqref{eq:func_delta} need to be estimated by their empirical counterparts.
 Denote by $F_{\ell j}^\star = \frac1{n_\ell} \sum_{k=1}^{n_\ell} \widehat \varepsilon_{\ell jk}$ and $G_j^\star = \frac1a \sum_{\ell = 1}^a F_{\ell j}^\star$ the wild bootstrap versions of $\widehat F_{\ell j} - F_{\ell j}$ and $\widehat G_{j} - G_{j}$, respectively.
 Finally, we obtain the following wild bootstrap counterpart of $\widehat{p}_{ij} - p_{ij}$:
 \begin{align}
 \label{eq:wild_bs_p}
  p_{ij}^\star = \int G_j^\star(x) \d \wh F_{ij}(x) - \int F_{\ell j}^\star(x) \d \widehat G_j(x) .
 \end{align}
 Combined into an $\R^{ad}$-vector $\b p^\star = \sqrt{N} (p_{11}^\star, p_{12}^\star, \dots, p_{ad}^\star)'$, 
 we have the following conditional central limit theorems which hold under both the null hypothesis $H_0^p(\b T) : \{\b T \b p = \b 0\}$ and the alternative hypothesis $H_1^p(\b T) : \{\b T \b p \neq \b 0\}$:
 }
 \begin{thm}
 \label{thm:cclt}
  Suppose Condition~\ref{cond:main} holds.
  As $N \rightarrow \infty$, we have, 
  conditionally on $\b X$, 
  \begin{align}
  \label{eq:bs_conv}
   \b p^\star \stackrel{d}{\longrightarrow} \b Z \sim N_{ad}( \bs 0_{ad}, \bs \Sigma)
  \end{align}
  in outer probability,
  where $\bs \Sigma$ is {\color{black} as in Theorem \ref{thm:clt}.}
 \end{thm}
 \begin{cor}
 \label{cor:ATS_wild}
  Suppose Condition~\ref{cond:main} holds.
  As $N \rightarrow \infty$, we have, 
  conditionally on $\b X$, 
  \begin{align}
  \label{eq:bs_conv_TN}
   T_N^\star = {\b p^\star}' \b T \b p^\star 
    \oDo  \b Z' \b T \b Z \gDg \sum_{{\g h}=1}^{ad} \nu_h Y_h^2
  \end{align}
  in outer probability, 
  i.e. the same limit distribution as in Corollary~\ref{cor:ATS}.
 \end{cor}
 {\bl
 The corollary again follows from the continuous mapping theorem applied to Theorem~\ref{thm:cclt}.
 The conditional central limit theorem~\eqref{eq:bs_conv_TN} is sufficient for providing random quantiles
 which converge in outer probability to the quantiles of the asymptotic distribution of the ATS under $H_0^p(\b T)$, i.e. $\sum_{{\g h}=1}^{ad} \nu_h Y_h^2$:
 repeated realizations of $T_N^\star$ are derived and their empirical quantiles serve as critical values for the hypothesis tests.
 
 Similarly, instead of a \emph{wild} bootstrap, a variant of the classical bootstrap \citep{efron1979bootstrap} may be applied to obtain a similar convergence result.
 The procedure is as follows:
 for each group $\ell=1,\dots,a$, we randomly draw $n_\ell$ $d$-dimensional data vectors from $\b X_{\ell 1}, \dots, \b X_{\ell n_\ell}$ with replacement to obtain \emph{bootstrap samples} $\b X_\ell^* = ((\b X_{\ell1}^*)', \dots, (\b X_{\ell n_\ell}^*)')'$.
 Denote their marginal empirical distribution functions as $F_{\ell j}^*, j=1,\dots,d$.
 Then, $F^*_{\ell j} - \widehat F_{\ell j}$  is the bootstrap counterpart of $\widehat F_{\ell j} - F_{\ell j}$ 
 and the respective bootstrap versions $\b p^*$ and $T_N^*$ of $\widehat {\b p}$ and $T_N$ are derived analogously to the wild bootstrap.
 Conditional central limit theorems analogous to~\eqref{eq:bs_conv} and~\eqref{eq:bs_conv_TN} hold too; see Appendix~\ref{sec:cbs} for details.
 We compare the performances of both proposed resampling procedures, i.e. the wild bootstrap and the classical bootstrap, in Section~\ref{sim} below.
 }

 %
 %
 %

{\bl
\section{Deduced Inference Procedures} \label{inf}

The previous considerations directly imply that consistent and asymptotic level $\alpha$ tests for $H_0^p(\b T):\{\b T \b p = \b 0\}$ are given by 
$$\varphi_N^\star = \1\{T_N > c^\star(\alpha)\} \text{ and } \varphi_N^* = \1\{T_N > c^*(\alpha)\},$$ where $c^*(\alpha)$ and $c^\star(\alpha)$ denote the $(1-\alpha)$ quantile of the classical bootstrap and wild bootstrap versions of $ T_N$ given $\b X$, i.e.\ of $T_N^\star$ in case of the wild bootstrap. Their finite sample performance  will be studied in Section~\ref{sim} below. 
As described in Section~\ref{mod}, these tests can be used to infer various global null hypotheses of interest about (nonparametric) main and interaction effects of interest which can straightforwardly be inverted to construct confidence regions for these nonparametric effects.

Moreover, the results derived in Section~\ref{sec:asy} also allow post-hoc analyses, i.e. 
subsequent multiple comparisons. 
To exemplify the typical paths of action we consider the one-way situation with $a$ independent groups and nonparametric effect size vectors $\b p_i=(p_{i1},\dots,p_{id})'$ in group $i, i=1,\dots,a$. 
If the global null hypothesis 
$$H_0^p(\b P_a \otimes \b I_d): \{\b p_1= \dots =\b p_a\}$$ of equal effect size vectors is rejected, one is usually interested in inferring 
\begin{itemize}
 \item[(i)] the (univariate) endpoints that caused the rejection, as well as
 \item[(ii)] the groups showing significant differences (all pairs comparisons).
\end{itemize}
The above questions directly translate to testing the univariate hypotheses
\bqan\label{eq:uni hyp}
H_{0j}^p: \{p_{1j} = \dots = p_{aj}\},\quad j=1,\dots, d
\eqan
in case of (i) and to an all pairs comparison given by multivariate hypotheses
\bqan\label{eq:all pairs}
H_{0i\ell}^p: \{\b p_i = \b p_\ell \}, \quad 1\leq i<\ell\leq a
\eqan
in case of (ii). Note that our derived methodology allows for testing these hypotheses in a unified way by performing tests on all univariate endpoints for \eqref{eq:uni hyp}
 and by selecting pairwise comparison contrast matrices for \eqref{eq:all pairs}. Therefore, a first naive approach would be to adjust the individual tests accordingly (e.g. by Bonferroni or Holm corrections) to ensure control of the family-wise error rate. 
However, note that the effect size vectors are defined via component-wise comparisons. This implies that the intersection of all $H_{0j}^p$ as well as the intersection of all $H_{0i\ell}^p$ is exactly given by the global null hypothesis $H_0^p(\b P_a \otimes \b I_d)$. 
Moreover, we can even test all subset intersections of $H_{0j}^p, j=1,\dots,d$ (or $H_{0i\ell}^p, 1\leq i< \ell\leq p$) by choosing adequate contrast matrices and performing the corresponding bootstrap procedures. Thus, both questions can even be treated (separately) by applying the closed testing principle of \cite{marcus1976closed}. 
This is a major advantage over existing inference procedures that are developed for testing null hypotheses formulated in terms of distribution functions \citep{ellis2017nonparametric}. In particular, since equality of marginals does not imply equality of multivariate distributions, 
the closed testing principle cannot be applied to the latter to answer question (i).

To ensure a reasonable computation time, the above approach is only applicable for small or moderate $p$ and $a$. However, some computation time can be saved by formulating a hierarchy on the questions (either for study-specific reasons or by weighing up the sizes of $a$ and $p$). For example, assume that (i) is more important than (ii). In this case we may start by applying the closed testing algorithm to test hypotheses $H_{0j}^p$ and subsequently only infer 
pair-wise comparisons on the significant univariate endpoints (instead of testing all multivariate $H_{0i\ell}^p$).
Contrary, assume that $d$ is much larger than $a$. Then it may be reasonable to first infer (ii) and subsequently consider (i) for the significant pairs.

}
\section{Simulations} \label{sim}

\subsection{Continuous Data}


For the one-way layout, data was generated similarly to the simulation study in \cite{Kon:2015}.
We considered $a=2$ treatment groups and $d \in \{4, 8\}$ endpoints as well as the following covariance settings:
\bqa
\text{Setting 1: }&& {\bf{V}}_1 = {\bf{I}}_d + 0.5 ({\bf{J}}_d-{\bf{I}}_d) = {\bf{V}}_2,  \\
\text{Setting 2: }&& \b V_1 = \left((0.6)^{|r-s|}\right)_{r, s =1}^d = \b V_2.  \\
\eqa
 Setting 1 represents a compound symmetry structure, while Setting 2 is an autoregressive covariance structure. 
Data was generated as
\bqa
\b X_{ik} = \b V_i^{1/2} \boldsymbol{\epsilon}_{ik}, ~i=1, 2;~ k=1, \dots, n_i,
\eqa
where $\b V_i^{1/2}$ denotes {\g a} square root of the matrix $\b V_i$, i.e., $\b V_i =\b V_i^{1/2} \cdot \b V_i^{1/2}$. 
The i.i.d.~random errors $\boldsymbol{\epsilon}_{ik}=(\epsilon_{ik1}, \dots, \epsilon_{ikd})'$ with mean $E(\boldsymbol{\epsilon}_{ik})= \vnull_d$ and $\Cov(\boldsymbol{\epsilon}_{ik}) = \b I_{d \times d}$ were generated by simulating independent standardized components $\epsilon_{iks}= (Y_{iks}- E(Y_{iks}))/(\Var(Y_{iks}))^{1/2}$
for various distributions of $Y_{iks}$. In particular, we simulated {\pu standard} normal {\g and} {\pu standard} lognormal distributed random variables.
We investigated balanced as well as unbalanced designs with sample size vectors $\b n^{(1)}=(10, 10)$, $\b n^{(2)} = (10, 20)$, and $\b n^{(3)}=(20, 10)$, {\g and increased sample sizes by adding $m \in \{0, 10, 30, 50\}$ to each element of the respective vector $\b n^{(h)}, h=1, 2, 3$.} 
In this setting, we tested the null hypothesis of no treatment effect $H_0^{p}: \{(\b P_a \kp \b I_d)\b p = \vnull_{2d} \}=\{\b p_1 = \b p_2\}$, where $\b p_i = (p_{i1}, \dots, p_{id})', i=1, 2$, {\pu and $\b p = (\b p_1', \b p_2')'$}.
{\g All simulations were conducted using the \textsc{R}-computing environment \citep{R}, version 3.2.3, each with 5,000 simulation runs and 5,000 bootstrap iterations.}

The results for {\g the} normal and lognormal distribution are displayed in Table \ref{tab:cont}. 
{\g The wild bootstrap approach shows a very good type-I error control for normally distributed data and $d=4$ dimensions, even for sample sizes as small as $\b n =(10, 10)'$. For the other scenarios considered, particularly for lognormal data, we need slightly larger sample sizes to achieve good type-I error control. 
An exception is the scenario with $d=8$ dimensions and covariance setting S2 with lognormal data, where the wild bootstrap approach maintains the pre-assigned level of 5\% already for sample sizes as small as 10.
	However, sample sizes of about 40 are enough to ensure very good type-I error rates across all scenarios considered here. }
The classical, group-wise bootstrap, in contrast, leads to slightly larger type-I error rates as compared to the wild bootstrap and only maintains the 5\% level for sample sizes of about 60.

\begin{table}[h]
	\centering
	\caption{Type-I error results for normal and lognormal distributed data with $d=4$ and $d=8$ dimensions, varying sample sizes and different covariance settings.}
		\label{tab:cont}
	\begin{tabular}{c|c|c|c|cccc|cccc}
		\hline
		& & & & \multicolumn{4}{c|}{wild bootstrap} &\multicolumn{4}{c}{group-wise bootstrap}\\
	&	distr & Cov & $\bs n$  &  $m=$ 0 & 10 & 30 & 50 &  0 &  10 & 30 & 50  \\ 
		\hline
	\multirow{12}{*}{$d=4$} &			\multirow{6}{*}{normal} & \multirow{3}{*}{S1} & (10, 10) &  5.2 & 5.5 & 5.5 & 5.5 & 8 & 5.8 & 5.8 & 5.5 \\ 
	& &  & (10, 20) & 5.1 & 5.1 & 5.4 & 5.2 & 7.7 & 6.5 & 5.4 & 4.9 \\ 
	& &  & (20, 10) & 5.5 & 5.3 & 4.9 & 5.3 & 7.1 & 6.4 & 5.2 & 5 \\ \cline{3-12}
	& & \multirow{3}{*}{S2} & (10, 10) & 4.8 & 5.2 & 5.2 & 5.3 & 7.8 & 5.8 & 5.9 & 5.2 \\ 
	& &  & (10, 20) & 5 & 4.9 & 5.6 & 4.9 & 7.7 & 6 & 5.7 & 5.1 \\ 
	& &  & (20, 10) & 5 & 5.3 & 5 & 5.1 & 6.9 & 6.3 & 5.2 & 4.9 \\ 	\cline{2-12}
	
	&	\multirow{6}{*}{lognormal} & \multirow{3}{*}{S1} & (10, 10) & 6 & 5.8 & 5.7 & 5.2 & 8.4 & 6.3 & 5.8 & 5.6 \\ 
	& &  & (10, 20) & 5.9 & 5.8 & 5.7 & 5.3 & 8.2 & 6.7 & 5.7 & 5.6 \\ 
	& &  & (20, 10) & 6.5 & 6.2 & 5.3 & 5.5 & 7.9 & 6.2 & 5.8 & 5.3 \\ \cline{3-12}
	& &  \multirow{3}{*}{S2} & (10, 10) & 5.8 & 6 & 5.7 & 5.1 & 8.2 & 6.6 & 5.9 & 5.3 \\ 
	& &  & (10, 20) & 5.5 & 5.7 & 5.5 & 5.3 & 8.2 & 6.4 & 5.6 & 5.2 \\ 
	& &  & (20, 10) & 6.2 & 6.1 & 4.9 & 5.3 & 7.8 & 6.4 & 5.9 & 5.3 \\ 
	\hline
		\multirow{12} {*}{$d=8$}				&	\multirow{6}{*}{normal} & \multirow{3}{*}{S1}  & (10, 10) & 6 & 5.3 & 5.2 & 5.1 & 8.1 & 6.3 & 5.7 & 6.1 \\ 
	& &  & (10, 20) & 5.7 & 5.9 & 5.2 & 4.7 & 7.7 & 5.6 & 5.4 & 5.6 \\ 
	& &  & (20, 10) & 5.9 & 5.3 & 4.8 & 4.8 & 7.5 & 5.7 & 5.8 & 5.7 \\ \cline{3-12}
	& &  \multirow{3}{*}{S2}& (10, 10) & 3.5 & 4.3 & 4.1 & 4.7 & 6.6 & 6.4 & 5.1 & 5.6 \\ 
	& &  & (10, 20) & 3.7 & 4.9 & 5 & 4.4 & 7 & 5.7 & 5.5 & 5.2 \\ 
	& &  & (20, 10) & 4.3 & 4.5 & 4.4 & 4.4 & 6.1 & 5.8 & 5.2 & 5.7 \\
	 \cline{2-12}
	&	\multirow{6}{*}{lognormal} & \multirow{3}{*}{S1} & (10, 10)& 6.9 & 5.9 & 5.3 & 5.7 & 8.6 & 6.5 & 6 & 5.6 \\ 
	& &  & (10, 20) & 6.4 & 6.7 & 5.7 & 4.6 & 8.4 & 6 & 6 & 5.5 \\ 
	& &  & (20, 10) & 6.5 & 6.2 & 5.1 & 4.7 & 8.3 & 6.1 & 5.9 & 5.4 \\ \cline{3-12}
	& &  \multirow{3}{*}{S2} & (10, 10) & 4.9 & 4.4 & 4.5 & 5.1 & 7.6 & 6.3 & 5.6 & 5.7 \\ 
	& &  & (10, 20) & 5.1 & 5.2 & 4.8 & 4.7 & 7.7 & 6.3 & 6 & 5.6 \\ 
	& &  & (20, 10) & 4.7 & 5 & 4.9 & 4.4 & 7.3 & 6.1 & 5.5 & 5.7 \\ 
				\hline
	\end{tabular}
\end{table}

\subsubsection{A heteroscedastic setting}

We simulated a heteroscedastic setting, where $H_0^p: \b T \b p = \vnull_{2d}$ is satisfied.
{\pu To this end, we took}
$$
{\bf{X}}_{ik} \sim N(\vnull_d, \sigma_i {\bf{I}}_d)
$$
for different choices of $\sigma_i \in \{1, 1.2, 2\}$ as well as sample sizes $n_i \in \{10, 20\}$. 
{\g Sample sizes were again increased as described above.}
The results are displayed in Table \ref{tab:growingSample}.
{\g In this case, we observe a rather conservative behavior across all scenarios, which improves with growing sample sizes, but is still slightly conservative in case of $d=8$ dimensions, even for sample sizes of 60 and 70. In this heteroscedastic setting, the classical, group-specific bootstrap yields better results in many scenarios, especially for growing dimension.
We note that this was the only studied setting where the group-specific bootstrap was better than the wild.
{\pu Apart from this exception,} it was the other way around.
}}

\begin{table}[h]
	\centering
	\caption{Type-I error in \% for the heteroscedastic setting.}
	\label{tab:growingSample}
	\begin{tabular}{c|c|c|cccc|cccc}
		\hline
			& & & \multicolumn{4}{c|}{wild bootstrap} &\multicolumn{4}{c}{group-wise bootstrap}\\
			&	$\sigma_i^2$ & $\bs n$  & $m=$0 & 10 & 30 & 50 &  0 &  10 & 30 & 50  \\ 
		\hline
	\multirow{9}{*}{$d=4$}			&	\multirow{3}{*}{(1, 2)} & (10, 10) & 1.3 & 2.4 & 3.5 & 4.4 & 6.2 & 5.2 & 5.3 & 5.6 \\ 
	& & (10, 20) & 1.4 & 2.5 & 3.5 & 3.6 & 5.6 & 5.8 & 5.2 & 4.9 \\ 
	& & (20, 10) & 2.2 & 3.9 & 3.9 & 4.5 & 6.6 & 6.1 & 6 & 5.2 \\ \cline{2-11}
		&	\multirow{3}{*}{(1, 1)} & (10, 10)  & 0.9 & 2.2 & 3.2 & 4.3 & 6.1 & 5.4 & 5.2 & 5.1 \\ 
	& & (10, 20) & 1.8 & 2.9 & 3.8 & 3.9 & 5.9 & 5.4 & 5.1 & 5.1 \\ 
	& & (20, 10) & 1.4 & 3.1 & 3.5 & 4 & 6 & 5.5 & 5.6 & 5.1 \\ \cline{2-11}
	&	\multirow{3}{*}{(1.2, 1)} & (10, 10) & 1 & 2.3 & 3.2 & 4.2 & 6.1 & 5.3 & 5.4 & 4.9 \\ 
	& & (10, 20) & 1.9 & 3.1 & 4 & 4 & 6.1 & 5.6 & 5.1 & 5 \\ 
	& & (20, 10) & 1.2 & 2.9 & 3.6 & 4 & 6 & 5.4 & 5.4 & 5.1 \\ 
		\hline
		\multirow{9}{*}{$d=8$}			&	\multirow{3}{*}{(1, 2)} & (10, 10)  & 0.3 & 1.5 & 2.8 & 3.7 & 4.7 & 4.7 & 4.6 & 4.8 \\ 
	& & (10, 20) & 0.3 & 1.7 & 3 & 3.7 & 4.2 & 4.8 & 4.9 & 5.2 \\ 
	& & (20, 10) & 1.3 & 2.8 & 4 & 4.4 & 4.4 & 4.9 & 4.8 & 4.8 \\ \cline{2-11}
	&	\multirow{3}{*}{(1, 1)} & (10, 10)  & 0.3 & 1.6 & 2.6 & 3.6 & 5 & 4.9 & 4.6 & 5.2 \\ 
	& & (10, 20) & 0.6 & 2.1 & 3.2 & 3.6 & 4.7 & 4.5 & 4.8 & 5.1 \\ 
	& & (20, 10) & 0.9 & 2.1 & 3.4 & 3.7 & 4.4 & 4.6 & 4.5 & 5 \\ \cline{2-11}
	&	\multirow{3}{*}{(1.2, 1)} & (10, 10)  & 0.2 & 1.5 & 2.8 & 3.5 & 4.8 & 4.9 & 4.8 & 5.3 \\ 
	& & (10, 20) & 0.8 & 2.2 & 3.3 & 3.8 & 4.8 & 4.3 & 4.7 & 5.1 \\ 
	& & (20, 10) & 0.8 & 2.1 & 3.3 & 3.6 & 4.4 & 4.6 & 4.5 & 4.9 \\ 
				\hline
	\end{tabular}
\end{table}

\subsection{Ordinal Data}

We simulated ordinal data using the function {\it ordsample} from the \textsc{R} package {\bf GenOrd} \citep{GenOrd,ferrari2012simulating}. 
{\g The package {\bf GenOrd} allows for simulation of discrete random variables with a given correlation structure and given marginal distributions.
	The latter are linked together via a Gaussian copula in order to achieve the desired correlation structure on the discrete components.
We simulated uniform marginal distributions, such that the outcomes in the $j$-th dimension are uniformly distributed on $j+1$ categories, $1 \leq j \leq d$. For the correlation structure, we used the same {\pu underlying} covariance matrices as in the continuous setting above.}
{\color{black} Again, w}e considered $d\in \{4, 8\}$ dimensions and the same sample sizes as above.
The results are displayed in Table \ref{tab:ordinal}, 
{\g showing a rather good type-I error control in the settings considered. The results here are similar to the ones obtained above for continuous data, with slightly larger type-I errors for the small sample scenarios.}
In this scenario, the group-wise bootstrap again shows a rather liberal behavior for small sample sizes. Although the behavior improves with growing sample size, the wild bootstrap is again superior here.

\begin{table}[h]
	\centering
	\caption{Type-I error rates in \% for ordinal data with different sample sizes and different covariance structures.}
	\label{tab:ordinal}
	\begin{tabular}{c|c|c|cccc|cccc}
		\hline
			& & & \multicolumn{4}{c|}{wild bootstrap} &\multicolumn{4}{c}{group-wise bootstrap}\\
			&	Cov.~setting & $\bs n$  &  $m=$0 & 10 & 30 & 50 &  0 &  10 & 30 & 50  \\ 
		\hline
			\multirow{6}{*}{$d=4$}	& \multirow{3}{*}{S1} & (10, 10) & 6.9 & 5.2 & 5.1 & 5.3 & 8.1 & 6.5 & 6.0 & 5.8 \\ 
		&	 & (10, 20) & 6.5 & 5.7 & 4.9 & 5.3 & 8.1 & 6.6 & 5.5 & 5.5 \\ 
			& & (20, 10) & 6.4 & 5.5 & 6.1 & 5.3 & 8.4 & 6.7 & 5.6 & 5.3 \\ \cline{2-11}
				& \multirow{3}{*}{S2} & (10, 10)  & 6.3 & 4.7 & 5.2 & 5.4 & 8.2 & 6.1 & 6.1 & 5.3 \\ 
				&& (10, 20) & 5.9 & 5.7 & 5.0 & 5.3 & 7.8 & 6.9 & 5.7 & 5.3 \\ 
				& & (20, 10) & 6.3 & 5.3 & 5.8 & 5.6 & 8.1 & 6.8 & 5.5 & 5.5 \\ 
				\hline
								
		\multirow{6}{*}{$d=8$}	& \multirow{3}{*}{S1} & (10, 10) & 6.4 & 5.8 & 5.1 & 4.8 & 8.7 & 6.7 & 5.9 & 5.6 \\ 
			& & (10, 20) & 6.4 & 5.7 & 5.5 & 5.6 & 7.5 & 7.0 & 6.2 & 6.0 \\ 
			& & (20, 10) & 6.4 & 5.5 & 5.2 & 5.0 & 8.0 & 6.5 & 5.3 & 6.0 \\ \cline{2-11}
			
			& \multirow{3}{*}{S2} & (10, 10) & 4.3 & 4.6 & 4.5 & 4.3 & 8.0 & 6.3 & 5.7 & 5.4 \\ 
			& & (10, 20) & 4.5 & 4.8 & 4.7 & 4.8 & 6.8 & 6.5 & 6.0 & 5.6 \\ 
			& & (20, 10) & 4.5 & 4.5 & 4.7 & 4.8 & 7.4 & 6.3 & 5.9 & 5.9 \\ 
				\hline
	\end{tabular}
\end{table}

In addition to type-I error rates, we have also compared the two bootstrap approaches with respect to their power behavior. The results can be found in Appendix~\ref{sec:power}.

\section{Data Example}\label{app}

As a data example, we consider the data set `marketing' in the \textsc{R}-package {\bf ElemStatLearn} \citep{ElemStatLearn}. This data set contains information on the annual household income along with 13 other demographic factors of shopping mall customers in the San Francisco Bay Area. Most of the variables in this data set are measured on an ordinal scale, {\pu rendering} mean-based approaches unfeasible. For our example, we consider the influence of sex 
 on annual household income and {\g educational status}. 
The annual household income is categorized in 9 categories ranging from `less than \$10,000' to `\$75,000 and more', while {\g education ranges} from `Grade 8 or less' to `Grad Study' {\g (6 categories)}.
This two-dimensional outcome is to be analyzed with respect to the influence factor sex (Male vs.~Female) 
.

The original data set consists of 8993 observations. After removing those observations with missing values in one of the variables considered here, {\g 8907} observations remain {\g (4041 male and 4866 female participants)}. 

{\g The estimated unweighted treatment effects are displayed in Table \ref{tab:effects}, while Figure \ref{fig:ecdf} shows the empirical distribution functions for the two dimensions for male and female participants, respectively. 
The effects can be interpreted in the following way: For both income and education, male participants tend to have higher values than female participants, i.e., males tend to have higher annual incomes and a higher level of education than females.}

A statistical analysis of the data example based on {\g our wild bootstrap approach} reveals {\g a} highly significant effects ($p${\pu -value} $< 0.0001$) of {\g sex on the two-dimensional outcome data, i.e., household income and educational status differ significantly between male and female participants.}

As a sensitivity analysis, we imputed the missing values using the \textsc{R} package \pkg{missForest} \citep{missForest, missForestArticle}. 
Performing the analysis as above with the imputed data leads to the same results.

\begin{table}[h]
	\centering
	\caption{\g Estimated treatment effects for the two dimensions Income and Education.}
	\label{tab:effects}
	\begin{tabular}{c|cc}
	Sex & Income & Education \\
	\hline
	Male & 0.511 & 0.517\\
	Female & 0.489 &  0.483 \\
	\end{tabular}
\end{table}

\begin{figure}[h]
	\includegraphics[width = \textwidth]{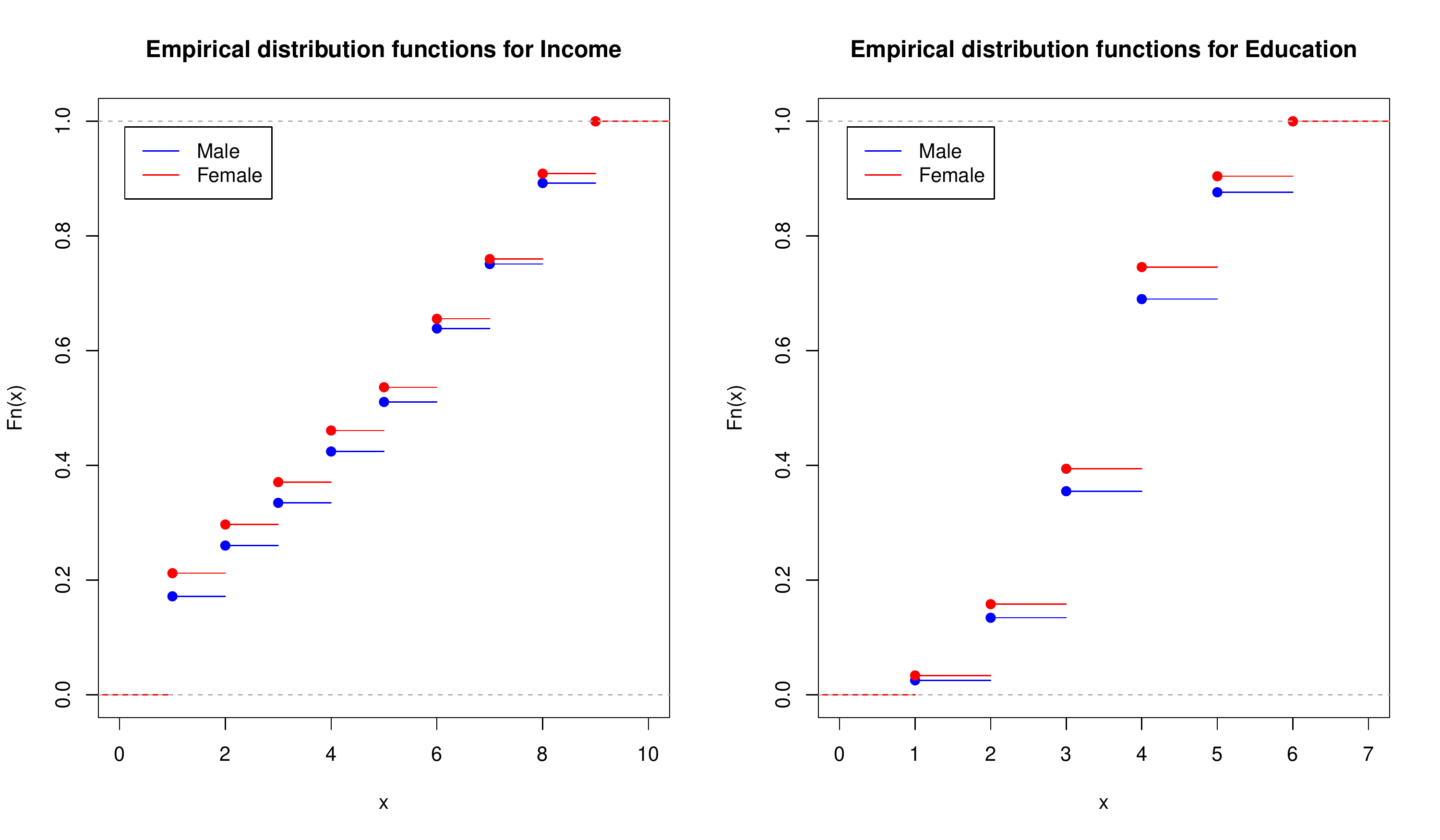}
	\caption{{\g Empirical distribution functions for male and female participants in the dimensions Income and Education, respectively.}
	}
	\label{fig:ecdf}
\end{figure}

%

\section{Conclusions and Discussion} \label{dis}

We have considered an extension of the unweighted treatment effects recently proposed by \citet{brunner2016JRSSB} to multivariate data. These effects do not depend on the sample sizes and allow for transitive ordering. We have rigorously proven the asymptotic behavior of the vector of unweighted treatment effects $\widehat{\b p}$ and proposed two bootstrap approaches to derive data-driven critical values for global and multiple test decisions. 
We proved the asymptotic validity of the bootstrap using empirical process arguments and analyzed its behavior in a large simulation study, where we considered continuous and ordinal distributions with different covariance settings and sample sizes. The {\g wild} bootstrap method performed {\g very} well in most scenarios {\g for small to moderate sample sizes. 
Only in heteroscedastic settings, it turned out to be rather conservative for small sample sizes. In this scenario, the sample-specific bootstrap provided better results. 
 
 {\g  In order to make them easily available for users, the proposed methods have been implemented {\pu by Sarah Friedrich}
	 in an \textsc{R} package {\bf rankMANOVA}, which is available from GitHub (https://github.com/smn74/rankMANOVA). 
 
{\pu In future work we will consider extensions of the present set-up
to censored multivariate data as well as address the question
``Which resampling method remains valid and performs preferably?''.
Here, a challenge will be the correct treatment of ties:
The wild bootstrap ceases to reproduce the correct limit distribution in case of right-censored and tied data
if it is not adjusted accordingly  \citep{dobler17}.
On the other hand, \cite{akritas86} has verified that Efron's bootstrap for right-censored data \citep{efron81} still works in the presence of ties.
The planned future paper may also be considered an extension of the article by \cite{dopa17} to the multi-sample and multivariate case.}

\section*{Acknowledgment} 
This work was supported by the German Research Foundation.


\appendix

\section{Proofs} \label{app:proofs}
{\bl Throughout, let $\mb P_{1, n_1}, \dots, \mb P_{a, n_a}$ be the empirical processes based on the samples $\b X_1, \dots, \b X_a$, respectively,
which are indexed by the class of functions $\mac G = \mac F \circ \Pi$,
where 
$$\mac F = \{ \1_{(-\infty, x]}(\cdot), \1_{(-\infty, x)}(\cdot) : x \in \R \}, $$
and $\Pi = \{ \pi_j : j=1,\dots, d \}$ is the class of all canonical coordinate projections $\pi_j : \R^d \rightarrow \R, (x_1, \dots, x_d) \mapsto x_j$.
Using this indexation, it is easily possible to derive the normalized empirical distribution functions $\widehat F_{ij}$ from $\mb P_{i, n_i}$.
We also see that every group-specific empirical process $\mb P_{i, n_i}$ can be considered as an element of $\ell^\infty(\mac G)$ which contains all bounded sequences with indices in $\mac G$.

}

\begin{proof}[Proof of Theorem~\ref{thm:clt}]
  Clearly, $\widehat{\b p}$ can be written as the image of all group- and component-specific, {\pu normalized} empirical distribution functions 
  $\frac1{n_i} \sum_{k=1}^{n_i} {\pu c(x-X_{ijk})}$, $j=1, \dots, d$, $i=1, \dots, a$, under a Hadamard-differentiable mapping {\pu $\phi$}.
  This can be seen following the lines in the proof of Theorem~2.1 in \cite{dopa17}
  where the case $a=2$ is discussed; in their proof, $K=\infty$ needs to be chosen.
  Hence, asymptotic normality follows from an application of the functional delta-method; cf. Theorem~3.9.4 in \cite{vaart96}.
  The asymptotic covariance structure of the resulting multivariate normal distribution is derived in detail in Appendix~\ref{app:cov}, where the asymptotic linear expansion of $\widehat{\b p}$ in all empirical distribution functions is utilized.
 \end{proof}

 \begin{proof}[Proof of Theorem~\ref{thm:cclt}]
 {\bl
 First note that, given $\b X$, we have conditional convergence in distribution 
  of $\b F_N^\star = \sqrt{N} (F_{11}^\star, F_{12}^\star, \dots, F_{ad}^\star)'$ to a multivariate Brownian bridge process in outer probability:
 this follows from an application of the conditional Donsker Theorem~3.6.13 in \cite{vaart96} in combination with Example~3.6.12 dedicated to the wild bootstrap and the choice of the Vapnik-\u{C}ervonenkis subgraph class
 $\mac F$ concatenated with the class of all canonical coordinate projections $\Pi$.
 This conserves the Vapnik-\u{C}ervonenkis subgraph property as argued in  Lemmata~2.6.17(iii) and 2.6.18(vii) of \cite{vaart96}.

 Next, recall the asymptotic linear representation~\eqref{eq:func_delta} of $\widehat{p}_{ij} - {p}_{ij}$
 which followed from the functional delta-method and which motivated the wild bootstrap version \eqref{eq:wild_bs_p}.
 This presentation involves Hadamard-derivatives that depend on estimated quantities:
 $$ \phi_{ij; \wh F}' : (\ell^\infty(\mac G))^a  \rightarrow \R, \quad (\textrm{P}_1, \dots, \textrm{P}_a) \mapsto \int \Big( \frac1a \sum_{\ell = 1}^a \textrm{F}_{\ell j} \Big) \d \widehat F_{ij} - \int \textrm{F}_{i j} \d \Big( \frac1a \sum_{\ell = 1}^a \widehat F_{\ell j} \Big). $$
 Here each $\textrm{P}_\ell, \ell=1,\dots, a,$ is a distribution on $\R^d$, which may hence be considered as an element of $\ell^\infty(\mac G)$,
 with marginal normalized distribution functions $\textrm{F}_{\ell j}$.

 We apply the extended continuous mapping theorem (Theorem 1.11.1 in \citealp{vaart96}) to the (random) functional $\phi'_{\wh F} = (\phi_{11; \wh F}', \phi_{12; \wh F}', \dots, \phi_{ad; \wh F}') : \ell^\infty(\mac G) \rightarrow \R^{ad}$.
 }
 The extended continuous mapping theorem is applied for almost every realization of $\b X$ thanks to the subsequence principle (Lemma~1.9.2 in \citealp{vaart96}):
 convergence in outer probability is equivalent to outer almost sure convergence along subsequences.
 The actual requirement for an application of the extended continuous mapping theorem is satisfied as well: 
 note that
 $\phi'_{\wh F}$ basically consists of integral mappings of the form 
 $$ \psi : D(\R) \times BV_1(\R) \rightarrow \R , \quad (f,g) \mapsto \int f \d g $$ 
 where $D(\R)$ is the space of right- (or left-)continuous functions on $\R$ with existing left- (or right-)sided limits
 and $BV_1(\R)$ is the subspace of functions with total variation bounded by 1.
 Lemma~3.9.17 in \cite{vaart96} states that $\psi$ is Hadamard-differentiable, hence continuous.
 We conclude that for all sequences of functions $(f_n)_{n \in \N}$ and $(g_n)_{n \in \N}$, which converge to $f_0$ in $D(\R)$ and to $g_0$ in $BV_1(\R)$, respectively, the sequence of functionals $\psi_n : f \mapsto \int f \d g_n$  satisfies $\psi_n(f_n) \rightarrow \int f_0 \d g_0$ as $n \rightarrow \infty$.
 {\bl All in all, the extended continuous mapping theorem, combined with the conditional central limit theorem for the wild bootstrapped empirical distribution functions as stated at the beginning of this proof, 
 concludes the proof of the conditional convergence in distribution of $\b p^\star$.}
 \end{proof}

\section{Covariances} \label{app:cov}
In this appendix we derive the covariance matrix $\bs \Sigma$ of the multivariate limit normal distribution in Theorem~\ref{thm:clt}.
The exact representation may not be strictly necessary for the practical purposes in this paper because a covariance estimator is not required due to the wild bootstrap asymptotics as described in Theorem~\ref{cor:ATS_wild}.
But the covariances below will give some insights into the asymptotically independent components of $\sqrt{N}(\wh {\b p} - \b p)$ and what kind of studentization may be applied if one wishes to test sub-hypotheses.
Furthermore, it is important to see that the limit distribution is not degenerate.
  Therefore, let $\widehat{\b w} $ be the vector consisting of all
  $$ \wh w_{\ell i j} = \int \wh F_{\ell j} \d \wh F_{ij} = \frac1{n_\ell} \frac1{n_i} \sum_{k=1}^{n_i} \sum_{r=1}^{n_\ell} c(X_{ijk} - X_{\ell j r}) 
  ; \quad \ell, i = 1, \dots, a, \ j=1,\dots, d.  $$
  This estimator is consistent for the vector, say, $\b w$ consisting of the different $w_{\ell ij}$.
  As an intermediate result, we are interested in the asymptotic covariance matrix 
  of the $\sqrt{N} (\wh {\b w} - \b w)$ vector, 
  i.e. in the limits $\sigma_{\ell i j, \ell' i' j'}$ of
  $$ N \cdot cov \Big( \int \wh F_{\ell j} \d \wh F_{ij} , \int \wh F_{\ell' j'} \d \wh F_{i'j'} \Big).$$
  To this end, we consider an asymptotically linear development 
  which is due to the functional delta-method: 
  Let ${\bl \psi} : (f,g) \mapsto \int f \d g$ {\bl again} denote the Wilcoxon functional; cf. Section~3.9.4.1 in \cite{vaart96}.
  As $N \rightarrow \infty$ and $\lim n_i / N \rightarrow \lambda_i$, \ $\lim n_\ell / N \rightarrow \lambda_\ell$ (according to Condition~\ref{cond:main}),
  \begin{align*}
   \sqrt{N} & \Big( \int \wh F_{\ell j} \d \wh F_{ij} - \int F_{\ell j} \d F_{ij} \Big)
    = \sqrt{N} ( {\bl \psi}(\wh F_{\ell j}, \wh F_{ij}) - {\bl \psi}(F_{\ell j}, F_{ij})) \\
   & = \sqrt{N} {\bl \psi'_{(F_{\ell j}, F_{ij})} } (\wh F_{\ell j} - F_{\ell j}, \wh F_{ij} - F_{ij}) + o_p(1) \\
   & = \int \sqrt{N} ( \wh F_{\ell j} - F_{\ell j}) \d F_{ij}
    +	\int \sqrt{N} F_{\ell j} \d ( \wh F_{\ell j} - F_{\ell j}) + o_p(1) \\
   & = \sqrt{N} \Big[ - \int F_{ij} \d \wh F_{\ell j} + \int F_{\ell j} \d \wh F_{ij} + \int F_{ij}\d F_{\ell j} - \int F_{\ell j} \d F_{ij} \Big] + o_p(1) \\
   & = \sqrt{N} \Big[ - \frac{1}{n_\ell} \sum_{r=1}^{n_\ell} F_{ij}(X_{\ell jr}) 
      + \frac{1}{n_i} \sum_{k=1}^{n_i} F_{\ell j}(X_{ijk}) + \int F_{ij}\d F_{\ell j} - \int F_{\ell j} \d F_{ij} \Big] + o_p(1) .
  \end{align*}
  Thus, we know that $\sigma_{\ell i j, \ell' i' j'}$ is the limit of
  \begin{align*}
    N \cdot cov & \Big( - \frac{1}{n_\ell} \sum_{r=1}^{n_\ell} F_{ij}(X_{\ell jr}) 
      + \frac{1}{n_i} \sum_{k=1}^{n_i} F_{\ell j}(X_{ijk}), \\
     & \qquad - \frac{1}{n_{\ell'}} \sum_{r'=1}^{n_{\ell'}} F_{i'j'}(X_{\ell' j'r'}) 
      + \frac{1}{n_{i'}} \sum_{k'=1}^{n_{i'}} F_{\ell' j'}(X_{i'j'k'}) \Big) \\
      & = \delta_{\ell \ell'} \frac{N}{n_\ell} cov(F_{ij}(X_{\ell j 1}), F_{i'j'}(X_{\ell j' 1}) )
       - \delta_{\ell i'} \frac{N}{n_\ell} cov(F_{ij}(X_{\ell j 1}), F_{\ell'j'}(X_{\ell j' 1}) ) \\
       & + \delta_{i i'} \frac{N}{n_i} cov(F_{\ell j}(X_{i j 1}), F_{\ell'j'}(X_{i j' 1}) )
       - \delta_{i \ell'} \frac{N}{n_i} cov(F_{\ell j}(X_{i j 1}), F_{i' j'}(X_{i j' 1}) ),
  \end{align*}
  where $\delta_{i i'} = 1\{i = i'\}$ is Kronecker's delta.
  We continue by calculating any of the above covariances, but we need to distinguish between two cases: \\
  Equal coordinates $\underline{j = j'}$:
  \begin{align*}
    cov(F_{ij}(X_{\ell j 1}), F_{i'j}(X_{\ell j 1}) )
    & = \int F_{ij}(u) F_{i'j}(u) \d F_{\ell j}(u) - \int F_{ij}(u) \d F_{\ell j}(u) \int F_{i'j}(u) \d F_{\ell j}(u) \\
    & = \tau_{ii'\ell j} - w_{i\ell j} w_{i' \ell j}.
  \end{align*}
  Unequal coordinates  $\underline{j \neq j'}$: Denote by $F_{\ell jj'}$ the joint normalized distribution function of $X_{\ell j 1}$ and $X_{\ell j' 1}$.
  \begin{align*}
    cov(F_{ij}(X_{\ell j 1}), F_{i'j'}(X_{\ell j' 1}) )
    & = \int F_{ij}(u) F_{i'j'}(v) \d F_{\ell j j'}(u,v) - \int F_{ij}(u) \d F_{\ell j}(u) \int F_{i'j'}(u) \d F_{\ell j'}(u) \\
    & = \rho_{ii'\ell j j'} - w_{i\ell j} w_{i' \ell j'}.
  \end{align*}
  Recall that $ w_{iij} = \frac12$.
  To sum up, we have the following asymptotic covariances (symmetric cases not listed):
  \begin{align*}
   \begin{cases}
    0 & \{i, \ell \} \cap \{i', \ell'\} = \emptyset \text{ or } i=\ell=i'=\ell' \\
    \frac{N}{n_\ell} (\tau_{ii \ell j} - w_{i \ell j}^2) 
      + \frac{N}{n_i} (\tau_{\ell \ell i j} - w_{\ell i j}^2)  
	  & j = j', i=i' \neq \ell = \ell' \\
     - \frac{N}{n_i} (\tau_{i \ell' i j} - w_{i i j} w_{\ell' i j}) 
      + \frac{N}{n_i} (\tau_{i \ell' i j} - w_{i i j} w_{\ell' i j})  
	  & j = j', i=i' = \ell \neq \ell' \\
      \frac{N}{n_\ell} (\tau_{i \ell \ell j} - w_{i \ell j} w_{\ell \ell j}) 
      - \frac{N}{n_\ell} (\tau_{i \ell \ell j} - w_{i \ell j} w_{\ell \ell j})  
	  & j = j', i \neq i' = \ell = \ell' \\
	- \frac{N}{n_\ell} (\tau_{ii\ell j} - w_{i\ell j}^2) 
	- \frac{N}{n_i} (\tau_{\ell \ell i j} - w_{\ell i j}^2) & j = j', i = \ell' \neq i' = \ell \\
	\frac{N}{n_i} (\tau_{\ell \ell' i j} - w_{\ell i j} w_{\ell' i j}) & j = j', i = i' \neq \ell \neq \ell' \neq i\\
	- \frac{N}{n_i} ( \tau_{\ell i' i j} - w_{\ell i j} w_{i' i j} ) & j = j', i = \ell' \neq i' \neq \ell \neq i \\
	\frac{N}{n_\ell} (\rho_{ii \ell jj'} - w_{i \ell j} w_{i \ell j'}) 
      + \frac{N}{n_i} (\rho_{\ell \ell i jj'} - w_{\ell i j} w_{\ell i j'})  
	  & j \neq j', i=i' \neq \ell = \ell' \\
     - \frac{N}{n_i} (\rho_{i \ell' i jj'} - w_{i i j} w_{\ell' i j'}) 
      + \frac{N}{n_i} (\rho_{i \ell' i jj'} - w_{i i j} w_{\ell' i j'})  
	  & j \neq j', i=i' = \ell \neq \ell' \\
      \frac{N}{n_\ell} (\rho_{i \ell \ell jj'} - w_{i \ell j} w_{\ell \ell j'}) 
      - \frac{N}{n_\ell} (\rho_{i \ell \ell jj'} - w_{i \ell j} w_{\ell \ell j'})  
	  & j \neq j', i \neq i' = \ell = \ell' \\
	- \frac{N}{n_\ell} (\rho_{ii\ell jj'} - w_{i\ell j} w_{i\ell j'}) 
	- \frac{N}{n_i} (\rho_{\ell \ell i jj'} - w_{\ell i j} w_{\ell i j'}) & j \neq j', i = \ell' \neq i' = \ell \\
	\frac{N}{n_i} (\rho_{\ell \ell' i jj'} - w_{\ell i j} w_{\ell' i j'}) & j \neq j', i = i' \neq \ell \neq \ell' \neq i\\
	- \frac{N}{n_i} ( \rho_{\ell i' i jj'} - w_{\ell i j} w_{i' i j'} ) & j \neq j', i = \ell' \neq i' \neq \ell \neq i \\
   \end{cases}
   \\ =
   \begin{cases}
    0 & \{i, \ell \} \cap \{i', \ell'\} = \emptyset \text{ or } i=\ell=i'=\ell' \\
     & \text{or } i=i' = \ell \neq \ell' \text{ or } i \neq i' = \ell = \ell' \\
    \frac{N}{n_\ell} (\tau_{ii \ell j} - w_{i \ell j}^2) 
      + \frac{N}{n_i} (\tau_{\ell \ell i j} - w_{\ell i j}^2)  
	  & j = j', i=i' \neq \ell = \ell' \\
	- \frac{N}{n_\ell} (\tau_{ii\ell j} - w_{i\ell j}^2) 
	- \frac{N}{n_i} (\tau_{\ell \ell i j} - w_{\ell i j}^2) & j = j', i = \ell' \neq i' = \ell \\
	\frac{N}{n_i} (\tau_{\ell \ell' i j} - w_{\ell i j} w_{\ell' i j}) & j = j', i = i' \neq \ell \neq \ell' \neq i\\
	- \frac{N}{n_i} ( \tau_{\ell i' i j} - w_{\ell i j} w_{i' i j} ) & j = j', i = \ell' \neq i' \neq \ell \neq i \\
	\frac{N}{n_\ell} (\rho_{ii \ell jj'} - w_{i \ell j} w_{i \ell j'}) 
      + \frac{N}{n_i} (\rho_{\ell \ell i jj'} - w_{\ell i j} w_{\ell i j'})  
	  & j \neq j', i=i' \neq \ell = \ell' \\
	- \frac{N}{n_\ell} (\rho_{ii\ell jj'} - w_{i\ell j} w_{i\ell j'}) 
	- \frac{N}{n_i} (\rho_{\ell \ell i jj'} - w_{\ell i j} w_{\ell i j'}) & j \neq j', i = \ell' \neq i' = \ell \\
	\frac{N}{n_i} (\rho_{\ell \ell' i jj'} - w_{\ell i j} w_{\ell' i j'}) & j \neq j', i = i' \neq \ell \neq \ell' \neq i\\
	- \frac{N}{n_i} ( \rho_{\ell i' i jj'} - w_{\ell i j} w_{i' i j'} ) & j \neq j', i = \ell' \neq i' \neq \ell \neq i \\
   \end{cases}
  \end{align*}
  In order to present the above covariances in a more compact matrix notation, 
  we introduce the following matrices:
  Denote by $\bs 0_{p \times q} \in \R^{p \times q}$ the $(p \times q)$-matrix of zeros,
  by $\bs 0_{r} \in \R^{r}$ the $r$-dimensional column vector of zeros,
  by $\bs \tau_{ii'\ell \cdot} = diag(\tau_{ii'\ell 1}, \dots, \tau_{ii'\ell d}) \in \R^{d \times d}$ the $(d \times d)$-diagonal matrices of $\tau$'s,
  by 
  $$\bs \rho_{ii'\ell \cdot \cdot} 
  = \begin{pmatrix}
     0 & \rho_{ii'\ell 12} & \rho_{ii'\ell 13} & \dots & \rho_{ii'\ell 1 d} \\
     \rho_{ii'\ell 21} & 0 & \rho_{ii'\ell 23} & \dots & \rho_{ii'\ell 2 d} \\
     \vdots & \vdots & \vdots & \ddots & \vdots \\
     \rho_{ii'\ell (d-1) 1} & \rho_{ii'\ell (d-1) 2} & \rho_{ii'\ell (d-1) 3} & \dots & \rho_{ii'\ell (d-1) d} \\
     \rho_{ii'\ell d1} & \rho_{ii'\ell d2} & \rho_{ii' \ell d 3} & \dots & 0 \\
    \end{pmatrix}  \in \R^{d \times d} $$
  the $(d \times d)$-matrices of $\rho$'s with zeros along the diagonal entries,
  and {\pu the vector of treatment effects between groups $i$ and $i'$} by $\b w_{ii'\cdot} = (w_{ii'1}, w_{ii'2}, \dots, w_{ii'd})^T \in \R^d$.
  Recall that, in the whole $w$-vector, we first first the $\ell$-value, then the $i$-value,
  so that we first go through the component index $j$.
  With the above notation, we thus obtain the following first block of the covariance matrix in which $\ell = \ell' = 1$
  which, for general indices $\ell$ and $\ell'$, we denote by $\bs \Sigma_{\ell \ell'} \in \R^{da \times da}$:
  \begin{align*}
   \bs \Sigma_{11}
   & = \frac{N}{n_1} \begin{pmatrix}
      \bs 0_{d \times d} & \bs 0_{d \times d} & \bs 0_{d \times d} & \dots & \bs 0_{d \times d} \\
      \bs 0_{d \times d} & \bs \tau_{221\cdot} & \bs \tau_{231\cdot} & \dots & \bs \tau_{2a1\cdot} \\
      \bs 0_{d \times d} & \bs \tau_{321\cdot} & \bs \tau_{331\cdot} & \dots & \bs \tau_{3a1\cdot} \\
      \vdots & \vdots & \vdots & \ddots & \vdots \\
      \bs 0_{d \times d} & \bs \tau_{a21\cdot} & \bs \tau_{a31\cdot} & \dots & \bs \tau_{aa1\cdot}
     \end{pmatrix}
     + \begin{pmatrix}
      \bs 0_{d \times d} & \bs 0_{d \times d} & \bs 0_{d \times d} & \dots & \bs 0_{d \times d} \\
      \bs 0_{d \times d} & \frac{N}{n_2} \bs \tau_{112\cdot} & \bs 0_{d \times d} & \dots & \bs 0_{d \times d} \\
      \bs 0_{d \times d} & \bs 0_{d \times d} & \frac{N}{n_3} \bs \tau_{113\cdot} & \dots & \bs 0_{d \times d} \\
      \vdots & \vdots & \vdots & \ddots & \vdots \\
      \bs 0_{d \times d} & \bs 0_{d \times d} & \bs 0_{d \times d} & \dots & \frac{N}{n_a} \bs \tau_{11a\cdot}
       \end{pmatrix} \\
    & + \frac{N}{n_1} \begin{pmatrix}
      \bs 0_{d \times d} & \bs 0_{d \times d} & \bs 0_{d \times d} & \dots & \bs 0_{d \times d} \\
      \bs 0_{d \times d} & \bs \rho_{221\cdot \cdot} & \bs \rho_{231\cdot \cdot} & \dots & \bs \rho_{2a1\cdot \cdot} \\
      \bs 0_{d \times d} & \bs \rho_{321\cdot \cdot} & \bs \rho_{331\cdot \cdot} & \dots & \bs \rho_{3a1\cdot \cdot} \\
      \vdots & \vdots & \vdots & \ddots & \vdots \\
      \bs 0_{d \times d} & \bs \rho_{a21\cdot \cdot} & \bs \rho_{a31\cdot \cdot} & \dots & \bs \rho_{aa1\cdot \cdot}
     \end{pmatrix}
     + \begin{pmatrix}
      \bs 0_{d \times d} & \bs 0_{d \times d} & \bs 0_{d \times d} & \dots & \bs 0_{d \times d} \\
      \bs 0_{d \times d} & \frac{N}{n_2} \bs \rho_{112\cdot \cdot} & \bs 0_{d \times d} & \dots & \bs 0_{d \times d} \\
      \bs 0_{d \times d} & \bs 0_{d \times d} & \frac{N}{n_3} \bs \rho_{113\cdot \cdot} & \dots & \bs 0_{d \times d} \\
      \vdots & \vdots & \vdots & \ddots & \vdots \\
      \bs 0_{d \times d} & \bs 0_{d \times d} & \bs 0_{d \times d} & \dots & \frac{N}{n_a} \bs \rho_{11a\cdot \cdot}
       \end{pmatrix} \\
    & - \frac{N}{n_1} \begin{pmatrix}
         \bs 0_d \\ \b w_{21\cdot} \\ \b w_{31\cdot} \\ \vdots \\ \b w_{a1\cdot}
        \end{pmatrix}
        \begin{pmatrix}
         \bs 0_d \\ \b w_{21\cdot} \\ \b w_{31\cdot} \\ \vdots \\ \b w_{a1\cdot}
        \end{pmatrix}^T
      - \begin{pmatrix}
        \bs 0_{d \times d} & \bs 0_{d \times d} & \bs 0_{d \times d} & \dots & \bs 0_{d \times d} \\
        \bs 0_{d \times d} & \frac{N}{n_2} \b w_{12\cdot} \b w_{12\cdot}^T & \bs 0_{d \times d} & \dots & \bs 0_{d \times d} \\
        \bs 0_{d \times d} & \bs 0_{d \times d} & \frac{N}{n_3} \b w_{13\cdot} \b w_{13\cdot}^T & \dots & \bs 0_{d \times d} \\
        \vdots & \vdots & \vdots & \ddots & \vdots \\
        \bs 0_{d \times d} & \bs 0_{d \times d} & \bs 0_{d \times d} & \dots & \frac{N}{n_a} \b w_{1a\cdot} \b w_{1a\cdot}^T
        \end{pmatrix}
  \end{align*}
  Note that the other $\bs \Sigma_{\ell \ell}$-matrices have a similar structure but with the $\bs 0_{d \times d}$-matrices in the $\ell$th block row and block column and with all $1$'s replaced by $\ell$'s. In the same way,
  \begin{align*}
   \bs \Sigma_{21} 
   & = - \frac{N}{n_2} \begin{pmatrix}
    \bs 0_{d \times d} & \bs \tau_{112\cdot} & \bs 0_{d \times d} & \dots & \bs 0_{d \times d} \\
    \bs 0_{d \times d} & \bs 0_{d \times d} & \bs 0_{d \times d} & \dots & \bs 0_{d \times d} \\
    \bs 0_{d \times d} & \bs \tau_{132\cdot} & \bs 0_{d \times d} & \dots & \bs 0_{d \times d} \\
    \bs 0_{d \times d} & \bs \tau_{142\cdot} & \bs 0_{d \times d} & \dots & \bs 0_{d \times d} \\
    \vdots & \vdots & \vdots & \ddots & \vdots \\
    \bs 0_{d \times d} & \bs \tau_{1a2\cdot} & \bs 0_{d \times d} & \dots & \bs 0_{d \times d}
   \end{pmatrix}
  - \frac{N}{n_1} \begin{pmatrix}
    \bs 0_{d \times d} & \bs \tau_{221\cdot} & \bs \tau_{231\cdot} & \dots & \bs \tau_{2a1\cdot} \\
    \bs 0_{d \times d} & \bs 0_{d \times d} & \bs 0_{d \times d} & \dots & \bs 0_{d \times d} \\
    \bs 0_{d \times d} & \bs 0_{d \times d} & \bs 0_{d \times d} & \dots & \bs 0_{d \times d} \\
    \bs 0_{d \times d} & \bs 0_{d \times d} & \bs 0_{d \times d} & \dots & \bs 0_{d \times d} \\
    \vdots & \vdots & \vdots & \ddots & \vdots \\
    \bs 0_{d \times d} & \bs 0_{d \times d} & \bs 0_{d \times d} & \dots & \bs 0_{d \times d}
   \end{pmatrix} \\
  & - \frac{N}{n_2} \begin{pmatrix}
    \bs 0_{d \times d} & \bs \rho_{112\cdot \cdot} & \bs 0_{d \times d} & \dots & \bs 0_{d \times d} \\
    \bs 0_{d \times d} & \bs 0_{d \times d} & \bs 0_{d \times d} & \dots & \bs 0_{d \times d} \\
    \bs 0_{d \times d} & \bs \rho_{132\cdot \cdot} & \bs 0_{d \times d} & \dots & \bs 0_{d \times d} \\
    \bs 0_{d \times d} & \bs \rho_{142\cdot \cdot} & \bs 0_{d \times d} & \dots & \bs 0_{d \times d} \\
    \vdots & \vdots & \vdots & \ddots & \vdots \\
    \bs 0_{d \times d} & \bs \rho_{1a2\cdot \cdot} & \bs 0_{d \times d} & \dots & \bs 0_{d \times d}
   \end{pmatrix}
  - \frac{N}{n_1} \begin{pmatrix}
    \bs 0_{d \times d} & \bs \rho_{221\cdot \cdot} & \bs \rho_{231\cdot \cdot} & \dots & \bs \rho_{2a1\cdot \cdot} \\
    \bs 0_{d \times d} & \bs 0_{d \times d} & \bs 0_{d \times d} & \dots & \bs 0_{d \times d} \\
    \bs 0_{d \times d} & \bs 0_{d \times d} & \bs 0_{d \times d} & \dots & \bs 0_{d \times d} \\
    \bs 0_{d \times d} & \bs 0_{d \times d} & \bs 0_{d \times d} & \dots & \bs 0_{d \times d} \\
    \vdots & \vdots & \vdots & \ddots & \vdots \\
    \bs 0_{d \times d} & \bs 0_{d \times d} & \bs 0_{d \times d} & \dots & \bs 0_{d \times d}
   \end{pmatrix} \\
     & + \begin{pmatrix}
    \bs 0_{d \times d} & \bs 0_{d \times d} & \bs 0_{d \times d} & \bs 0_{d \times d} & \dots & \bs 0_{d \times d} \\
    \bs 0_{d \times d} & \bs 0_{d \times d} & \bs 0_{d \times d} & \bs 0_{d \times d} & \dots & \bs 0_{d \times d} \\
    \bs 0_{d \times d} & \bs 0_{d \times d} & \frac{N}{n_3} \bs \tau_{213\cdot} & \bs 0_{d \times d} & \dots & \bs 0_{d \times d} \\
    \bs 0_{d \times d} & \bs 0_{d \times d} & \bs 0_{d \times d} & \frac{N}{n_4} \bs \tau_{214\cdot} & \dots & \bs 0_{d \times d} \\
    \vdots & \vdots & \vdots & \vdots & \ddots & \vdots \\
    \bs 0_{d \times d} & \bs 0_{d \times d} & \bs 0_{d \times d} & \bs 0_{d \times d} & \dots & \frac{N}{n_a} \bs \tau_{21a\cdot}
   \end{pmatrix}
   \\
   & + \begin{pmatrix}
    \bs 0_{d \times d} & \bs 0_{d \times d} & \bs 0_{d \times d} & \bs 0_{d \times d} & \dots & \bs 0_{d \times d} \\
    \bs 0_{d \times d} & \bs 0_{d \times d} & \bs 0_{d \times d} & \bs 0_{d \times d} & \dots & \bs 0_{d \times d} \\
    \bs 0_{d \times d} & \bs 0_{d \times d} & \frac{N}{n_3} \bs \rho_{213\cdot \cdot} & \bs 0_{d \times d} & \dots & \bs 0_{d \times d} \\
    \bs 0_{d \times d} & \bs 0_{d \times d} & \bs 0_{d \times d} & \frac{N}{n_4} \bs \rho_{214\cdot \cdot} & \dots & \bs 0_{d \times d} \\
    \vdots & \vdots & \vdots & \vdots & \ddots & \vdots \\
    \bs 0_{d \times d} & \bs 0_{d \times d} & \bs 0_{d \times d} & \bs 0_{d \times d} & \dots & \frac{N}{n_a} \bs \rho_{21a\cdot \cdot}
   \end{pmatrix} \\
   & + \frac{N}{n_2} 
   \begin{pmatrix}
    \bs 0_{da \times d} & 
   \begin{pmatrix}
     \b w_{12\cdot}^T \\ \bs 0_d^T \\ \b w_{32\cdot}^T \\ \b w_{42\cdot}^T \\ \vdots \\ \b w_{a2\cdot}^T
   \end{pmatrix} \otimes \b w_{12\cdot}
    & \bs 0_{da \times d(a-2)}
   \end{pmatrix} \\
   & + \frac{N}{n_1} 
   \begin{pmatrix}
    \begin{pmatrix}
    \bs 0_{d}^T & \b w_{21\cdot}^T & \b w_{31\cdot}^T & \dots & \b w_{a1\cdot}^T 
    \end{pmatrix} \otimes \b w_{21\cdot}
    \\
    \bs 0_{d(a-1) \times da}
   \end{pmatrix}
   \\
   & - \begin{pmatrix}
      \bs 0_{d \times d} & \bs 0_{d \times d} & \bs 0_{d \times d} & \bs 0_{d \times d} & \dots & \bs 0_{d \times d} \\
      \bs 0_{d \times d} & \bs 0_{d \times d} & \bs 0_{d \times d} & \bs 0_{d \times d} & \dots & \bs 0_{d \times d} \\
      \bs 0_{d \times d} & \bs 0_{d \times d} & \frac{N}{n_3} \b w_{23\cdot} \b w_{13\cdot}^T & \bs 0_{d \times d} & \dots & \bs 0_{d \times d} \\
      \bs 0_{d \times d} & \bs 0_{d \times d} & \bs 0_{d \times d} & \frac{N}{n_4} \b w_{24\cdot} \b w_{14\cdot}^T & \dots & \bs 0_{d \times d} \\
      \vdots & \vdots & \vdots & \vdots & \ddots & \vdots \\
      \bs 0_{d \times d} & \bs 0_{d \times d} & \bs 0_{d \times d} & \bs 0_{d \times d} & \dots & \frac{N}{n_a} \b w_{2a\cdot} \b w_{1a\cdot}^T
     \end{pmatrix}
  \end{align*}
  The representation of the general block matrix $\bs \Sigma_{\ell \ell'}$ with $\ell \neq \ell'$ is similarly obtained,
  where the zero-rows have to shifted to the row block number $\ell$ and the zero-column to the column block number $\ell'$.
  Furthermore, the repeating $1$'s and $2$'s in the above representation need to be replaced with $\ell'$'s and $\ell$'s, respectively.

  Since each $\wh {p}_{ij}$ is the mean of $\wh w_{1ij}, \wh w_{2ij}, \dots, \wh w_{aij}$,
  we conclude that the limit covariance matrix of $\sqrt{N} ( \wh {\b p} - \b p)$ is given by
  \begin{align*}
  \bs \Sigma = \lim_{N \rightarrow \infty} \frac{1}{a^2} \sum_{\ell = 1}^a \sum_{\ell' = 1}^a \bs \Sigma_{\ell \ell'}. 
  \end{align*}
  \qed


\section{Theory for the classical, group-wise bootstrap}
\label{sec:cbs}

{\g 
In this section, we present the results of the classical, group-wise bootstrap \pu{as a competitor of the wild bootstrap}.
The bootstrap is implemented in the following fashion: For each group $i =1,\dots , a$, we draw $n_i$ independent selections $\b X_{ik}^* = (X_{i1k}^*, \dots,  X_{idk}^*)'$
randomly with replacement from the vectors $\b X_{ik} , k = 1,\dots , n_i$. These are used to build the bootstrapped empirical
distribution functions $F_{ij}^*$ as described in Section 2, i.e., $F_{ij}^* (x) = \frac{1}{n_i}\sum_{k=1}^{n_i}c(x-X_{ijk}^*)$.
 The bootstrapped treatment effects are then obtained via
$$
p_{ij}^* = \int G_j^* dF_{ij}^* = \frac{1}{a}\sum_{\ell = 1}^a w_{\ell ij}^*.
$$
\pu{The asymptotics of Theorem \ref{thm:clt} and Corollary~\ref{cor:ATS}} for these bootstrapped treatment effects 
\pu{hold under both the null hypothesis and the alternative hypothesis
which is shown by} using similar arguments as in the proofs for the wild bootstrap:

\begin{thm}
 \label{thm:cclt_classic}
  Suppose Condition~\ref{cond:main} holds.
  As $N \rightarrow \infty$, we have,
  conditionally on $\b X$, 
  \begin{align*}
   \sqrt{N}(\b p^* - \widehat{\b p}) \stackrel{d}{\rightarrow} \b Z \sim N_{ad}( \bs 0_{ad}, \bs \Sigma)
  \end{align*}
  in outer probability,
  where $\bs \Sigma$ is {\color{black} as in Theorem \ref{thm:clt}.}
 \end{thm}
 The continuous mapping theorem immediately implies the corresponding conditional convergence in distribution for the bootstrapped ATS:
 \begin{cor}
  Suppose Condition~\ref{cond:main} holds.
  As $N \rightarrow \infty$, we have, conditionally on $\b X$, 
  under both,  $H_0^p$ and ($\b T \b p \neq \b 0$),
  \begin{align*}
   T_N^* = N ( {\b p}^* - \widehat{\b p} )' \b T ( {\b p}^* - \widehat{\b p} ) 
    \oDo  \b Z' \b T \b Z \gDg \sum_{{\g h}=1}^{ad} \nu_h Y_h^2
  \end{align*}
  in outer probability, 
  i.e. the same limit distribution as in Corollary~\ref{cor:ATS}.
 \end{cor}

 \begin{proof}[Proof of Theorem~\ref{thm:cclt_classic}]
  Similarly, as argued in the proof of Theorem~\ref{thm:clt}, $\b p^*$ is obtained 
  as a Hadamard-differentiable functional of all bootstrapped (normalized) empirical distribution functions 
  $F_{ij}^*(t) = \frac1{n_i} \sum_{k=1}^{n_i} c(t - X^*_{ijk})$, $j=1, \dots, d$, $i=1, \dots, a$.
  As the conditional central limit theorem holds in outer probability for each bootstrapped empirical distribution function, i.e. for each
  $$ \sqrt{n_i} (F_{ij}^*(t) - \wh F_{ij}(t)) = \frac{1}{\sqrt{n_i}} \Big( \sum_{k=1}^{n_i} c(t - X^*_{ijk}) - \sum_{k=1}^{n_i} c(t - X_{ijk}) \Big),$$
  cf. Theorem~3.6.1 in \cite{vaart96},
  the convergence is transferred to $\sqrt{N}(\b p^* - \widehat{\b p})$ by means of the functional delta-method for the bootstrap; cf. Theorem~3.9.11 in \cite{vaart96}.
 \end{proof}

Note that the pooled bootstrap, corresponding to drawing with replacement from the combined sample,
 is not available in this context:
 this method involves asymptotic eigenvalues other than those in Corollary~\ref{cor:ATS_wild}.
 Hence, the limiting distribution of the pooled bootstrapped ATS would only be correct in special cases.

	\section{Additional simulation results: Power}
	\label{sec:power}
	
	In order to in compare the power behavior of the two bootstrap methods, we have considered a shift alternative, i.e., we simulated data as
	\bqa
	\tilde{\b X}_{ik} = \bs \mu_i + \b X_{ik}, ~i=1, 2;~ k=1, \dots, n_i,
	\eqa
	where $\bs \mu_1 = \b 0_{d}$ and $\bs \mu_2 = (\delta, \dots, \delta)'$ for $\delta \in \{0, 0.5, 1, 1.5, 2, 3\}$ and $\b X_{ik}$ corresponds to the respective random vectors simulated in Section \ref{sim} of the paper.
	The simulation results for some scenarios are exemplarily shown in Figures \ref{fig:power1} and \ref{fig:power2}. We find that the power results for both bootstrap approaches are almost identical in the chosen situations.
	
	\begin{figure}[h]
		\includegraphics[width =\textwidth]{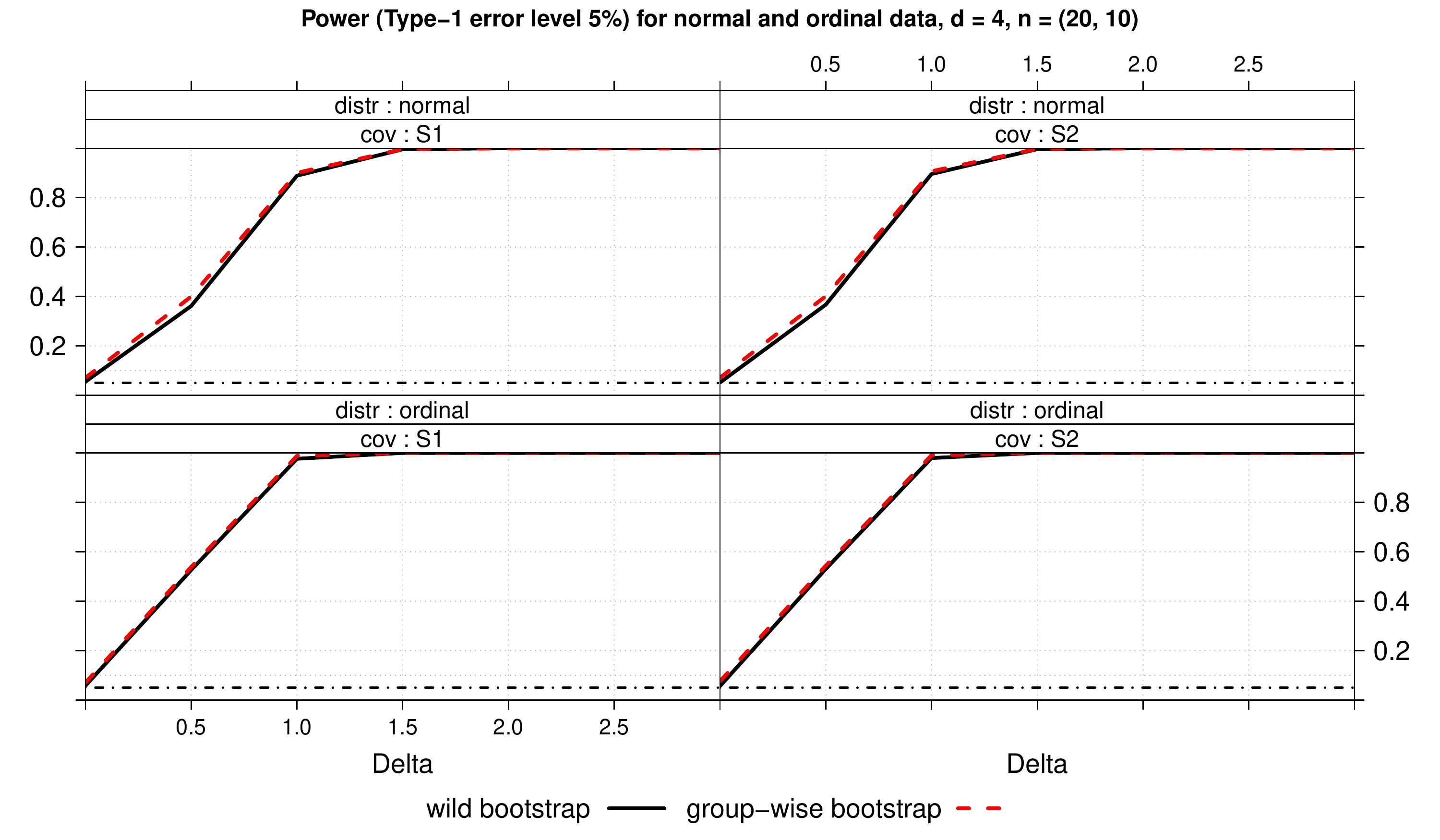}
		\caption{Power simulation results for continuous (normally distributed) and ordinal data with $d=4$ dimensions and $\b n=(20, 10)'$.}
		\label{fig:power1}
	\end{figure}

		\begin{figure}[h]
			\includegraphics[width =\textwidth]{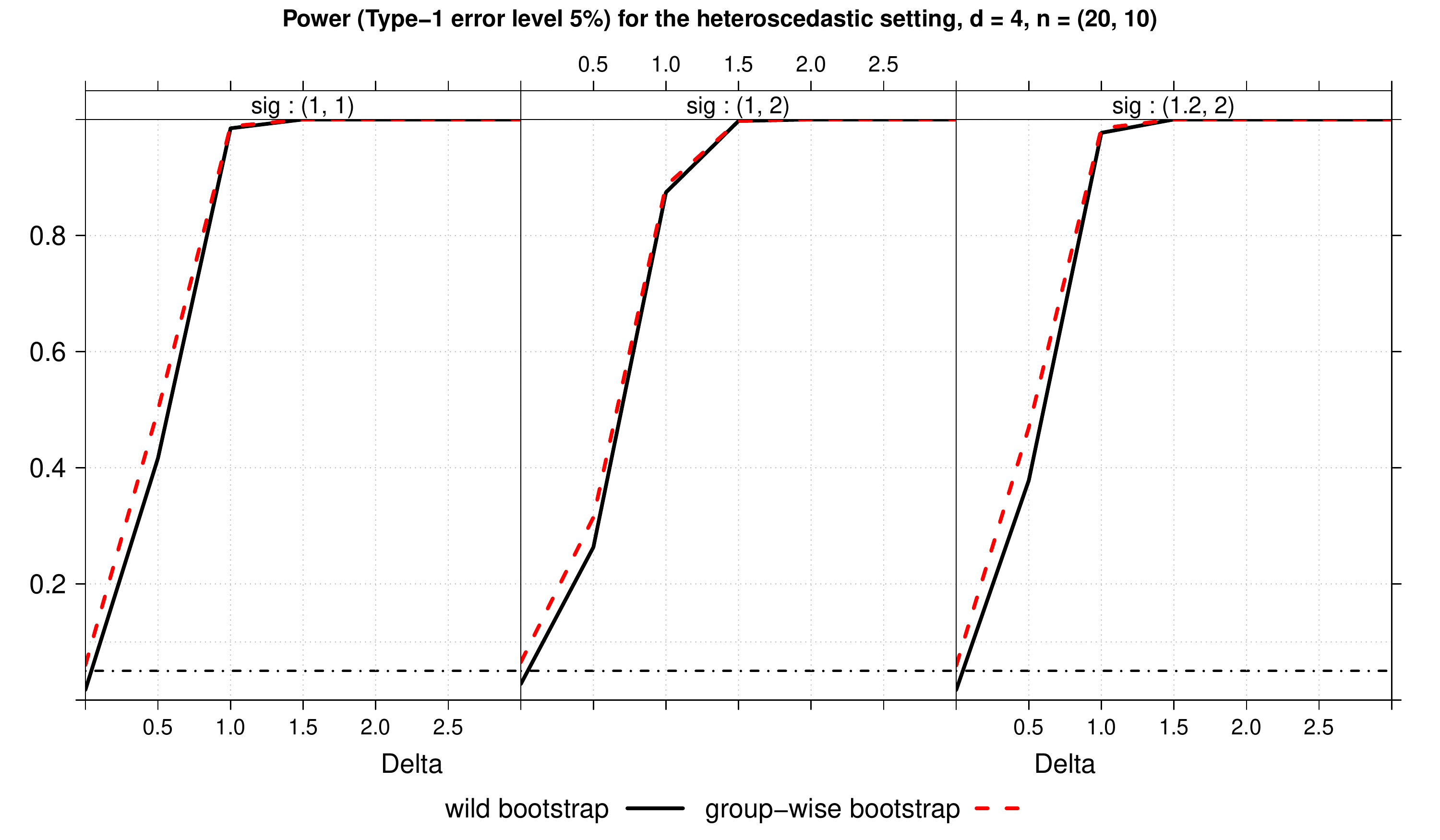}
			\caption{Power simulation results for heteroscedastic data with $d=4$ dimensions and $\b n=(20, 10)'$.}
			\label{fig:power2}
		\end{figure}

\section{Additional analyses of the data example}
\label{sec:add_ana}
In order to demonstrate the proposed methods in a more complex context, we consider another analysis of the 'marketing' data example. Additionally to the variables described in Section \ref{app} we now also include the factor 'language', which describes the language spoken most often at home with the three levels 'English', 'Spanish' and 'Other'. After removing the observations with missing values in these variables, 8423 observations remain. We thus want to conduct some exploratory analyses of the two-dimensional outcome in a two-way layout with factors 'sex' and 'language'.
The estimated nonparametric effects are displayed in Table \ref{table:trteffects} and the wild bootstrap approach results in highly significant $p$-values for the two main as well as the interaction effect, see Table \ref{tab:pval}.

Since the interaction hypothesis is significant, we continue by analyzing male and female participants separately.
In order to further interpret the results, we also apply the post hoc comparisons described in Section \ref{inf}. 
In particular, since the global null hypothesis is significant in both groups, we continue with the pairwise comparisons of the languages. Since again all results are significant at 5 \% level, we finally consider the univariate outcomes. The results are displayed in Table \ref{table:pairwise}.
This reveals some interesting aspects of the data. For example, the significant difference between `English' and `Other' in the male group is driven by education, while there is no significant effect on the income. A similar result is obtained for `Spanish' vs.~`Other' in the female group.

\begin{table}[ht]
	\centering
	\caption{\g Estimated nonparametric effects for the two dimensions Income and Education.}
	\label{table:trteffects}
	\begin{tabular}{cc|cc}
		Sex & Language & Income & Education \\
		\hline
		\multirow{3}{*}{Male} & English& 0.586 & 0.605\\
		& Spanish & 0.560 &0.568\\
		& Other & 0.464& 0.407\\ \hline
		\multirow{3}{*}{Female} & English & 0.401 &0.367\\
		& Spanish & 0.529 &0.554\\
		& Other & 0.460& 0.499\\		
	\end{tabular}
\end{table}

\begin{table}[ht]
	\centering\caption{Multivariate and univariate $p$-values for the main and interaction effects.}
\label{tab:pval}
	\begin{tabular}{c|ccc}
		Effect	& \multicolumn{3}{c}{$p$-value}\\
	 & multivariate & Income & Education \\\hline
		Sex & $<0.0001$ & $<0.0001$& $<0.0001$\\
	Language & $<0.0001$ & $<0.0001$& $<0.0001$\\
					Sex:Language & $<0.0001$& $<0.0001$& $<0.0001$\\
	\end{tabular}
\end{table}

\begin{table}[ht]
	\centering
	\caption{\g Pairwise comparisons of the different groups with respect to the two-dimensional outcome.}
	\label{table:pairwise}
	\begin{tabular}{cc|ccc}
			Sex	& Language & \multicolumn{3}{c}{$p$-value}\\
		&	& multivariate & Income & Education \\\hline
		\multirow{4}{*}{Male} & global null hypothesis & $<0.0001$ & $<0.0001$& $<0.0001$\\\cline{2-5}
		& English vs.~Spanish & $<0.0001$ & $<0.0001$& $<0.0001$\\
		& English vs.~Other & 0.028 &  0.067 &  0.029\\
		& Spanish vs.~Other & $<0.0001$ & 0.042 & $<0.0001$\\
 \hline
		\multirow{3}{*}{Female} & global null hypothesis & $<0.0001$ &$<0.0001$ & $<0.0001$\\ \cline{2-5}
		& English vs.~Spanish &$<0.0001$ & $<0.0001$&$<0.0001$\\
		& English vs.~Other &$<0.0001$  & $<0.0001$ &0.018\\
		& Spanish vs.~Other & $<0.0001$& 0.069& $<0.0001$\\		
	\end{tabular}
\end{table}

All these analyses can be conducted with the \textsc{R} package {\bf rankMANOVA} by splitting the data accordingly. The implementation of a routine for these calculations is part of future research.

 \newpage

\bibliography{Literatur}{}
\bibliographystyle{apalike}

\end{document}